\documentclass[11pt,oneside,dvipsnames]{amsart}
\usepackage{stackrel}
\usepackage[french, english]{babel}
\usepackage{a4wide}
\usepackage[utf8,latin1]{inputenc}
\usepackage{times}
\usepackage{adjustbox}
\usepackage{float}
\usepackage[cyr]{aeguill}
\usepackage{amsmath,amscd}
\usepackage{amssymb}
\usepackage{mathrsfs}
\usepackage{amsthm}
\usepackage{geometry}
\usepackage{graphicx}
\usepackage{epstopdf}
\usepackage{amsfonts}
\usepackage{amsmath}
\usepackage{amsmath,mathtools}
\usepackage{amsmath,graphicx}
\usepackage{amscd}
\usepackage{dsfont}
\usepackage{latexsym}
\usepackage{verbatim}
\usepackage{comment}
\usepackage{tikz}
\usepackage{tikz-cd}
\usetikzlibrary{matrix,arrows}
\usetikzlibrary{tikzmark}
\usepackage{tensor}
\usepackage{leftidx}
\usepackage[retainorgcmds]{IEEEtrantools}
\usepackage[title]{appendix}
\usepackage[mathscr]{euscript}
\usepackage[colorinlistoftodos]{todonotes}
\usepackage{bbm}
\usepackage{url}
\usepackage{upgreek}
\usepackage{mathtools, nccmath}
\usepackage[utf8]{inputenc}
\usepackage{multirow}
\usepackage{rotating}
\usepackage{marginnote}
\usepackage{hyperref} 
\usepackage{nameref}
\usepackage{ytableau}
\usepackage{youngtab}

\usepackage[all]{xy}
\usepackage[T1]{fontenc}

\usepackage{hyperref}
\usepackage{array}
\usepackage{tabularx}
\usepackage{yfonts}
\usepackage{setspace}  
 \usepackage[foot]{amsaddr}
\makeatletter
\renewcommand\subsubsection{\@startsection{subsubsection}{3}{\z@}%
                                     {-2ex\@plus -1ex \@minus -.2ex}%
                                     {-0.5em}
                                     {\normalfont\normalsize\bfseries}}
                                     
                                     \renewcommand\subsection{\@startsection{subsection}{3}{\z@}%
                                     {-3.25ex\@plus -1ex \@minus -.2ex}%
                                     {-0.5em}
                                     {\normalfont\normalsize\bfseries}}

\makeatother

\theoremstyle{plain}
\newtheorem{Thm}{Theorem}
\newtheorem{ThmA}{Theorem}

\newtheorem{Prop}[Thm]{Proposition}
\newtheorem{Lem}[Thm]{Lemma}
\newtheorem{Cor}[Thm]{Corollary}

\theoremstyle{definition}

\newtheorem{Rem}[Thm]{Remark}

\newtheorem{Defn}[Thm]{Definition}

\newtheoremstyle{named}{}{}{}{}{\bfseries}{.}{.5em}{\thmnote{#3}#1}
\theoremstyle{named}

\DeclareMathOperator{\Ad}{Ad}

\DeclareMathOperator{\Ch}{Ch}

\DeclareMathOperator{\gl}{\mathfrak{gl}}
\DeclareMathOperator{\GL}{GL}

\DeclareMathOperator{\Hom}{Hom}

\DeclareMathOperator{\Id}{Id}
\DeclareMathOperator{\Ima}{Im}

\DeclareMathOperator{\lisom}{\!\!\smash{\begin{array}{c}\sim\\[-1em]
\longrightarrow\end{array}}\!\!}

\DeclareMathOperator{\Lie}{Lie}

\DeclareMathOperator{\PGL}{PGL}

\DeclareMathOperator{\Rep}{Rep}

\DeclareMathOperator{\SL}{SL}
\DeclareMathOperator{\SO}{SO}
\DeclareMathOperator{\Sp}{Sp}

\DeclareMathOperator{\Spin}{Spin}

\DeclareMathOperator{\Stab}{Stab}

\DeclareMathOperator{\ds}{\!/\mkern-2mu/\mkern-2mu}

\newcommand*\circled[1]{
  \tikz[baseline=(char.base)]\node[shape=circle,draw,inner sep=0.2pt,font=\tiny,minimum size=8pt] (char) {#1};}

\newcommand{\hooklongrightarrow}{\lhook\joinrel\longrightarrow}

\makeatletter
\newcommand*{\rom}[1]{\expandafter\@slowromancap\romannumeral #1@}
\makeatother

\renewcommand{\descriptionlabel}[1]{\hspace\labelsep\upshape\bfseries #1.}
\makeatletter
\let\orgdescriptionlabel\descriptionlabel
\renewcommand*{\descriptionlabel}[1]{%
  \let\orglabel\label
  \let\label\@gobble
  \phantomsection
  \edef\@currentlabel{#1}%
  \let\label\orglabel
  \orgdescriptionlabel{#1}%
}
\makeatother

\usepackage{subfig}	 

\title{The most singular loci of twisted character varieties}
\author{Cheng Shu}
\address[C. Shu]{School of Mathematical Sciences, Tongji University, Shanghai, China}
\email[C. Shu]{shuc@tongji.edu.cn}

\usepackage{pxfonts}
\usepackage{indentfirst}
\usepackage{setspace}  
\usepackage{hyperref}
\colorlet{ivory}{Apricot!30!}
\colorlet{space}{black!85!}
\definecolor{bgc}{RGB}{29, 44, 46}
\definecolor{txt}{RGB}{223, 222, 189}
\definecolor{cmd}{RGB}{206, 151, 88}
\begin{document}
\let\bs\boldsymbol
\begin{abstract}
For any simply connected almost simple group $G$ and any compact Riemann surface of genus $g>1$, we study the singular locus of the character variety which parametrises fundamental group representations with a twist by a central element $z$ of $G$. While for $z=1$ the most singular representations have values in a maximal torus (and are invariant for the Weyl group action) and for generic $z$ every $z$-twisted representation is irreducible, for a general $z$ we find that the most singular representations have values in a particular Levi subgroup $L_z$ up to conjugation. From the explicit Lie-theoretic description of $L_z$, we determine the structure of the maximal stabiliser groups, and prove that the most singular loci is a disjoint union of generic type $A$ character varieties. 
\end{abstract}
\maketitle

\tableofcontents
\addtocontents{toc}{\protect\setcounter{tocdepth}{-1}}
\setcounter{tocdepth}{1}
\numberwithin{Thm}{section}
\numberwithin{equation}{section}
\addtocontents{toc}{\protect\setcounter{tocdepth}{1}}
\section{Introduction}\label{INT}

The moduli space of $G$-local systems on a compact Riemann surface, where $G$ is a reductive group over $\mathbb{C}$, is a Poisson variety and carries a natural symplectic structure on its smooth locus according to a classical result of Goldman \cite{G1}. There have been recent attempts to clarify the nature of the singularities of these affine varieties, and it is now understood that character varieties have symplectic singularities in the sense of Beauville \cite{Beau}; see \cite{BS} and \cite{Shu4}. In short, this means that the structure of the singularities, which is recorded in the resolutions of singularities, is compatible with the symplectic structure on the smooth locus in a suitable manner. As a corollary, character varieties are rational Gorenstein. Moreover, analyses of the codimensions of the singular loci show that character varieties have terminal singularities in many cases.

A more subtle question is whether character varieties are $\mathbb{Q}$-factorial or even locally factorial, which is crucial for the (non)existence of symplectic resolutions. The answer is rooted in how factoriality behaves under a GIT quotient. Indeed, character varieties are defined as the GIT quotients of some usually locally factorial varieties, called the representation varieties, and the question is really about whether local factoriality descends. The standard technique for such a problem is a theorem of Dr\'ezet \cite{Dre1}, later reformulated by Bellamy-Schedler in \cite{BS0} (see also \cite{Shu4}). This theorem says that, under reasonable conditions, if for every $G$-equivariant line bundle $\mathcal{L}$ on the representation variety and every point $x$ with a closed orbit, the stabiliser $\Stab_G(x)$ acts trivially on the fibre $\mathcal{L}_x$, then local factoriality descends. 

It turns out that there are stabilisers that act trivially on any equivariant line bundles for an obvious reason. Namely, if a stabiliser is a semi-simple group, then it has no nontrivial character (homomorphism to the multiplicative group $\mathbb{G}_m$), and thus acts trivially on the fibre of any line bundle. This is the case, when $G$ is semi-simple and $x$ is the point representing the trivial local system, so that $\Stab_G(x)=G$. Then, a method introduced in \cite{KLS} can be used to extend the triviality of the stabiliser action away from $x$. 

However, if $G$ is a semi-simple group that is not simply connected, then the associated character variety is not connected, and the nonneutral connected components do not contain the trivial local system as a point. This motivates the problem that we address in this article (up to isogeny): what is the maximal stabiliser in a nonneutral connected component of a representation variety? The results of this article will be used in a future version of \cite{Shu4} to complete the proof of the $\mathbb{Q}$-factoriality of character varieties.

To simply the formulation, we will work with the simply connected cover of $G$ and consider a twisted version of character variety. In the rest of this introduction, we assume that $G$ is a simply connected semi-simple group and $z\in Z_G$ denotes a central element. Let $\Pi$ be the fundamental group of a compact Riemann surface of genus $g\ge 2$. Define the associated representation variety:
$$
\Rep_z(\Pi,G):=\{(A_i,B_i)_{1\le i\le g}\in G^{2g}\mid \prod_{i=1}^g[A_i,B_i]=z\},
$$
where $[A_i,B_i]=A_iB_iA_i^{-1}B_i^{-1}$ is the commutator. If $z=1$, then this variety parametrises homomorphisms $\Pi\rightarrow G$, which we call $G$-valued representations in this article; if $z\ne 1$, then elements of $\Rep_z(\Pi,G)$ will be called $z$-twisted representations. By abuse of notaion, a $z$-twisted representation is also denoted by $\rho:\Pi\rightarrow G$. The affine GIT quotient
$$
\Ch_z(\Pi,G):=\Rep_z(\Pi,G)\ds G
$$
will be called the $z$-twisted character variety. We will write $\Rep(\Pi,G)$ and $\Ch(\Pi,G)$ if $z=1$.

\begin{ThmA}\label{Thm-A}
If $G$ is not of type $C_{2m+1}$, then there is a linearly reductive subgroup $G_z\subset G$ satisfying the following conditions:
\begin{itemize}
\item there exists $\rho\in \Rep_z(\Pi,G)$ with a closed $G$-orbit such that $\Stab_G(\rho)=G_z$; and
\item for any $\rho\in \Rep_z(\Pi,G)$ with a closed $G$-orbit, if $\Stab_G(\rho)$ contains $ G_z$, then $\Stab_G(\rho)=G_z$.
\end{itemize}  
If $G\cong\Sp_{2(2m+1)}$, then there is a linearly reductive subgroup $G'_z\subset G$ satisfying the same conditions with $G_z$ replaced by $G'_z$. Moreover, the identity connected components of $G_z$ and $G_z'$ have the largest possible rank among stabilisers and are semi-simple.
\end{ThmA}

The fact that these stabiliser groups have semi-simple identity components is what we need to prove the $\mathbb{Q}$-factoriality of character varieties. The proof of this theorem is a combination of Lemma \ref{Lem-Gz-Lz}, Corollary \ref{Cor-Stab=G_z} and Corollary \ref{Cor-Stab=G_z'}. The precise isomorphism type of the identity component of $G_z$ and $G'_z$ for every pair $(G,z)$ will be given by Theorem \ref{Thm-Cent-H_z} and Proposition \ref{Prop-C2m+1-Stab-circ}, and the component group of $G_z$ and $G'_z$ will be given by Theorem \ref{Thm-Gamma-G-z} and Proposition \ref{Prop-Gamma-G'-z}.

It is easy to see that if $z=1$, then the maximal stabiliser group is simply $G$, and the elements of $\Rep(\Pi,G)$ having $G$ as its stabiliser are precisely the $Z_G$-valued representations (i.e., elements of $Z_G^{2g}$, or $Z_G$-local systems). Therefore, the locus of representations with the maximal stabiliser is a set of isolated points, and its connected components form a finite abelian group. To extend such a statement to the case $z\ne 1$, we need the notion of Luna stratification. For any affine variety $X$ with an action of a reductive group $G$, the Luna stratification is a finite partition of $X\ds G$ into locally closed subvarieties:
$$
X\ds G=\bigsqcup_{(H)}(X\ds G)_{(H)}
$$
where $H\subset G $ is a linearly reductive subgroup and $(H)$ denotes the set of the $G$-conjugates of $H$. Every point $x\in (X\ds G)_{(H)}$ corresponds to a closed orbit in $X$, and the stabiliser of a point in this orbit is conjugate to $H$. Beware that in general this partition may not satisfy the axiom of stratification; that is, the closure of a stratum is a union of strata. In addition, each part $(X\ds G)_{(H)}$ may neither be connected nor smooth.

For quiver varieties, Bellamy-Schedler proved that the above partition is indeed a stratification and each stratum $(X\ds G)_{(H)}$ is connected and smooth; see \cite[Theorem 3.14 and Proposition 3.15]{BS0}. More interesting is the fact that this stratification is precisely the stratification by symplectic leaves; in particular, every stratum carries an induced symplectic structure.

For character varieties, however, we will prove the following.
\begin{ThmA}(Theorem \ref{Thm-sing-loci=Ch-A})\label{Thm-B}
If $G$ is not of type $C_{2m+1}$, then there is a semi-simple subgroup $H_z\subset G$ of type $A$, which contains $z$ as a generic central element, such that the connected components of $\Ch_z(\Pi,G)_{(G_z)}$ are in bijection with $(Z_{G_z}/Z_{G_z}\cap H_z)^{2g}$, and every connected component is isomorphic to $\Ch_z(\Pi,H_z)$. If $G\cong\Sp_{2(2m+1)}$, then there is a semi-simple subgroup $H'_z\subset G$ of type $A$, which contains $z$ as a generic central element, such that the connected components of $\Ch_z(\Pi,G)_{(G'_z)}$ are in bijection with $Z_{G_z}^{2g}$, and every connected component is isomorphic to $\Ch_z(\Pi,H'_z)$.
\end{ThmA}

This theorem says that the most singular locus (in the sense of Luna stratification) of a twisted character variety for a simply connected group $G$ is a disjoint union of generic type $A$ character varieties, and the connected components form a finite abelian group. For the convenience of the reader, we list $H_z$ as follows.
\begin{center}
\begin{tabular}{|c|c|c|c|c|c|c|c|c|c|}
\hline
$G$ & $A_{\ell}$ & $B_{\ell}$ & $C_{2m}$ & $C_{2m+1}$ & $D_{\ell}$ & $D_{2m+1}$ & $D_{2m}$ & $E_6$ & $E_7$ \\
\hline
$z$ & order $d$ & $\ne 1$ & $\ne 1$ & $\ne1$ & $\alpha_1$ & $\alpha_{\ell-1}$ or $\alpha_{\ell}$ & $\alpha_{\ell-1}$ or $\alpha_{\ell}$ & $\ne1$ & $\ne1$ \\ 
\hline
$H_z$ & $\SL_d$ & $\SL_2$ & $\SL_2$ & $\SL_2\times\SL_2$ & $\SL_2$ & $\SL_2\times\SL_4$ & $\SL_2$ & $\SL_3$ & $\SL_2$ \\
\hline
\end{tabular}
\end{center}
For $G\cong\Sp_{2(2m+1)}$, we have $H'_z\cong\SL_2$ which is diagonally embedded in $H_z$.

We briefly explain the strategy for the proofs of these two theorems. For a $z$-twisted representation $\rho=(A_i,B_i)_{1\le i\le g}$, let $H$ be the Zariski closure of the subgroup of $G$ generated by $\{A_i,B_i\}_i$. If the $G$-orbit of $\rho$ is closed, then $H$ is a linearly reductive subgroup $H\subset G$. Such a group is always an irreducible subgroup of some Levi subgroup, so that the latter is minimal among Levi subgroups containing $H$. The fact that $\rho$ is a $z$-twisted representation is an additional constraint on the possible Levi subgroups. Indeed, the derived subgroup of a relevant Levi subgroup must contain $z$ by the definition of a $z$-twisted representation. For example, a maximal torus cannot contain $H$ if $z\ne 1$. It turns out that for every $z\in Z_G$, there is a Levi subgroup $L_z$, unique up to conjugation, that is minimal among Levi subgroups satisfying $z\in[L,L]$. The Dynkin diagram of $L_z$ will be identified with a subdiagram of the Dynkin diagram of $G$; details are given in \S \ref{sec-BFM-Levi}.

It is perhaps a coincidence that every connected component of the Dynkin diagram of $L_z$ is of type $A$, so that $[L_z,L_z]$ is in fact a direct product of special linear groups (Levi subgroups of simply connected groups are simply connected). For $G=\SL_n$ and $z\in Z_G$ a root of unity of order $d$, we will see that $[L_z,L_z]$ is a direct product of $m$ copies of $\SL_d$, where $md=n$. Moreover, every direct factor of $z\in Z_{[L_z,L_z]}$ is a primitive $d$-th root of unity. In this situation, we say that $z$ is a generic central element, which implies that no proper Levi subgroup contains $z$ and, in particular, the $z$-twisted $L_z$-character variety is smooth. An obvious way to produce representations with large stabilisers is to require a twisted representation $\rho$ to have values in the diagonal subgroup of $[L_z,L_z]$, i.e., the $\SL_d$ diagonally embedded in $[L_z,L_z]$. In other words, we are taking a direct sum of $m$ copies of a $d$-dimensional irreducible representation. Classical theory says that the identity component of the automorphism group of such a representation is $\SL_m$. Note that scalar matrices $\zeta\Id$ with $\zeta^n=1$ but $\zeta^m\ne1$ gives rise to nontrivial connected components of the automorphism/stabiliser group. Our definition of $H_z$ or $H'_z$ is a natural generalisation of this picture to other types of $G$, and $G_z$ or $G_z'$ will be the centraliser of $H_z$ or $H'_z$ respectively. By computing the invariant vectors of $H_z$ or $H'_z$ acting on the Lie algebra of $G$, we determine the root system of the identity connected component of $G_z$ or $G'_z$. With a bit more effort, we will be able to show that these are the largest possible stabiliser groups.

Interestingly, the Levi subgroups $L_z$ already appeared in Borel-Friedman-Morgan's study of character varieties for $g=1$, playing a different role though. If $z=1$ and $g=1$, then $\Ch(\Pi,G)$ parametrises commuting semi-simple elements in $G$ up to conjugation. Indeed, there is an isomorphism of schemes (see \cite{LNY})
$$
\Ch(\Pi,G)\cong (T\times T)/W,
$$
where $T\subset G$ is a maximal torus and $W$ is the Weyl group of $G$ defined by $T$ acting diagonally on $T\times T$. Note that $T=L_z$ if $z=1$. If $z\ne 1$, then $\Ch_z(\Pi,G)$ parametrises almost commuting semi-simple elements, i.e., semi-simple elements $(s,t)\in G\times G$ with $[s,t]=z$. Again, a minimal Levi subgroup containing $s$ and $t$ is conjugate to $L_z$.  If $S_z$ is the connected centre of $L_z$ and $(s_0,t_0)\in L_z\times L_z$ is a given almost commuting pair, then we obtain other almost commuting pairs by multiplication by $S_z$. In fact, all almost commuting pairs are obtained this way; see \cite[Proposition 6.2.1]{LNY}. By introducing additional equivalence relations on $S_z\times S_z$, we have the following identification:
$$
\Ch_z(\Pi,G)\cong (\bar{S}_z\times\bar{S}_z)/W_z,
$$
where $\bar{S}_z=S_z/S_z\cap[L_z,L_z]$ and $W_z$ is a certain Weyl group acting diagonally on $\bar{S}_z\times\bar{S}_z$ (in most cases, it is indeed the Weyl group of $G_z$ in Theorem \ref{Thm-A}). It is a nontrivial theorem of Li-Nadler-Yun that this identification is in fact an isomorphism of schemes. To summarise, the character varieties for $g=1$ are completely determined by the Levi subgroups $L_z$. The results of the present article will reveal the role of $L_z$ in describing the most singular loci of the character varieties for $g>1$.

\subsection*{AI disclosure}\hfill

We have used ChatGPT-6 Astra Pro to test or verify some Lie-theoretic computations. We take full responsibility for the final content of this work.

\subsection*{Acknowledgement}\hfill

This work is supported by the China Postdoctoral Science Foundation under Grant Numbers 2025M773090 and 2025T180845 that are associated to the author's position at Westlake University, Hangzhou.

\section{Preliminaries}\label{sec-Pre}
\subsection{Reductive groups}\label{subsec-alg-grp}\hfill

For any affine algebraic group $G$, we denote by $G^{\circ}$ its identity connected component. A reductive group is understood to be connected. Let $G$ be a reductive group. There is an isogeny $Z_G^{\circ}\times [G,G]\rightarrow G$, where $[G,G]$ is the derived subgroup of $G$ and $Z_G$ is the centre of $G$, and the map is defined by the multiplication of $G$. It follows that there is a natural surjective homomorphism
\begin{equation}
\mathbf{D}_G:G\longrightarrow Z_G^{\circ}/Z_G^{\circ}\cap [G,G].
\end{equation}
Note that $Z_G^{\circ}\cap [G,G]$ is finite. When $G=\GL_n$, the map $\mathbf{D}_G$ is precisely the determinant map $\GL_n\rightarrow \mathbb{C}^{\ast}$.

Let $T\subset G$ be a maximal torus. The character lattice (resp. the cocharacter lattice) defined by $T$ is the free abelian group $X(T):=\Hom(T,\mathbb{G}_m)$ (resp. $Y(T):=\Hom(\mathbb{G}_m,T)$). The set of roots of $G$ with respect to $T$ will be denoted by $\Phi_G$, and the set of coroots by $\Phi_G^{\vee}$. The root lattice $Q_G$ (resp. coroot lattice $Q_G^{\vee}$) is the lattice spanned by $\Phi_G$ (resp. $\Phi_G^{\vee}$). Choose a basis $\Delta=\{\alpha_i\}_{1\le i\le\ell}$ of $\Phi_G$ and denote by $\{\varpi_i^{\vee}\}_{1\le i\le\ell}\subset\mathbb{Q}\otimes_{\mathbb{Z}}Y(T)$ the dual basis consisting of fundamental coweights, and by $\{\varpi_i\}_{1\le i\le \ell}\subset\mathbb{Q}\otimes_{\mathbb{Z}}X(T)$ the basis consisting of fundamental weights that are dual to the set of simple coroots $\Delta^{\vee}:=\{\alpha_i^{\vee}\}_{1\le i\le\ell}$. The weight lattice $P_G$ (coweight lattice $P_G^{\vee}$) is the lattice spanned by $\{\varpi_i\}_{1\le i\le \ell}$ (resp. $\{\varpi_i^{\vee}\}_{1\le i\le\ell}$). The subscript $G$ will be omitted if no confusion arises.

Suppose that $G$ is semi-simple. We will denote by $G_{sc}$ the simply-connected cover of $G$ and by $G_{ad}$ the adjoint quotient of $G$. If $T\subset G$ is a maximal torus, then its inverse image $T_{sc}\subset G_{sc}$ is a maximal torus and its image $T_{ad}\subset G_{ad}$ is also a maximal torus. 

\subsection{Central elements and the Dynkin diagram}\label{subsec-FOSS}\hfill 

We fix an injective group homomorphism $\iota:(\mathbb{Q}/\mathbb{Z})\hookrightarrow \mathbb{C}^{\ast}$. It should be thought of as the exponential map and its image consists of the torsion elements in $\mathbb{C}^{\ast}$. We will denote by $\tilde{\iota}$ the composition of $\iota$ and the natural projection $\mathbb{Q}\rightarrow\mathbb{Q}/\mathbb{Z}$. If $T$ is a torus, then we define $\tilde{\iota}_T:\mathbb{Q}\otimes_{\mathbb{Z}}Y(T)\rightarrow T(\mathbb{C})$ as sending $r\otimes\lambda$ to $\lambda(\tilde{\iota}(r))$. The image of $\tilde{\iota}_T$ consists of the finite order elements of $T$. The restriction of $\tilde{\iota}_T$ to $Y(T_{ad})$ induces an isomorphism $Y(T_{ad})/Y(T)\cong Z_G$. Now, suppose that $G$ is simply connected, so that $Y(T)=Q^{\vee}$. Since $Y(T_{ad})=P^{\vee}$, each fundamental coweight $\varpi_i^{\vee}$ gives a central element $\tilde{\iota}_T(\varpi_i^{\vee})$ of $G$. Let $\beta$ be the highest root for the basis $\Delta$, then $\beta=\sum_in_i\alpha_i$ for some $n_i\in\mathbb{Z}_{>0}$. The set $\{\varpi_i^{\vee}\mid n_i=1\}$ is a set of representatives of $Y(T_{ad})/Y(T)$. Therefore, for each nontrivial $z\in Z_G$, there is a unique $\varpi_i^{\vee}$ with $n_i=1$ corresponding to $z$. We will frequently represent elements of $Z_G$ by elements of $\Delta$ under this correspondence.

\subsection{Irreducible and completely reducible subgroups}\hfill

Let $G$ be a reductive group and let $H\subset G$ be a closed subgroup. We say that $H$ is an irreducible subgroup of $G$ if there is no proper parabolic subgroup of $G$ containing $H$. We say that $H$ is a completely reducible subgroup of $G$ if for any parabolic subgroup $P$ containing $H$, there exists a Levi factor $L$ of $P$ containing $H$. An irreducible subgroup of $G$ is completely reducible. Recall that a linear algebraic group is linearly reductive if every rational representation of it is semi-simple. In characteristic zero, the following are equivalent for a closed subgroup $H\subset G$:
\begin{itemize}
\item $H^{\circ}$ is reductive.
\item $H$ is linearly reductive.
\item $H$ is a completely reducible subgroup of $G$. 
\end{itemize}
See, for example, \cite[Section 2.2 and Section 2.4]{BMR}. 

It is well-known that in characteristic zero, the centraliser $C_G(H)$ is linearly reductive if $H$ is so.
\begin{Lem}\label{Lem-S<C_G(H)}
The following assertions hold:
\begin{itemize}
\item[(i)] Let $L\subset G$ be a Levi subgroup and let $H\subset L$ be an irreducible subgroup. Then, $S:=Z_L^{\circ}$ is a maximal torus of $C_G(H)^{\circ}$. 
\item[(ii)] Let $H\subset G$ be a linearly reductive subgroup and let $S\subset C_G(H)^{\circ}$ be a maximal torus. Then, $H$ is an irreducible subgroup of the Levi subgroup $L:=C_G(S)$.
\end{itemize}
\end{Lem}
\begin{proof}
Obviously, $S$ is a subtorus of $C_G(H)^{\circ}$. Let $T$ be a maximal torus of $C_G(H)^{\circ}$ containing $S$. Then, $H$ is contained in the Levi subgroup $M=C_G(T)$. If $T\ne S$, then $M$ is a proper Levi subgroup of $L$, which contradicts the assumption that $H$ is an irreducible subgroup of $L$. To prove (ii), suppose otherwise that there exists a proper parabolic subgroup $Q$ of $L$ that contains $H$. This implies that $H$ is contained in a Levi factor $M$ of $Q$. Since $M$ is a proper Levi subgroup of $L$, its connected centre $Z_M^{\circ}$ is strictly larger than $S$, contradicting the assumption that $S\subset C_G(H)^{\circ}$ is a maximal torus.
\end{proof}

\subsection{Irreducible and completely reducible representations}\label{subsec-rep-char}\hfill

The representation variety $\Rep_z(\Pi,G)$ parametrises $z$-twisted $G$-representations $\rho:\Pi\rightarrow G$. Such a $\rho$ is represented by a tuple $(A_i,B_i)_{1\le i\le g}\in G^{2g}$. The image of $\rho$, denoted by $\Ima\rho$, is defined to be the subgroup of $G$ generated by the set $\{A_i,B_i\mid 1\le i\le g\}$. The conjugation action of $G$ on itself induces a diagonal action of $G$ on $G^{2g}$ preserving the subvariety $\Rep_z(\Pi,G)$. The affine GIT quotient $\Ch_z(\Pi,G)$ parametrises the closed orbits in $\Rep_z(\Pi,G)$. Let $\rho\in \Rep_z(\Pi,G)$ and denote by $H$ the Zariski closure of $\Ima\rho$ and by $\Stab_G(\rho)$ the stabiliser group of $\rho$. The orbit of $\rho$ is closed if and only if $H$ is a completely reducible subgroup of $G$, and the orbit of $\rho$ is closed with $\Stab_G(\rho)^{\circ}=Z_G^{\circ}$ if and only if $H$ is an irreducible subgroup of $G$; in these cases, we will call $\rho$ a completely reducible or irreducible representation accordingly. The standard reference for these facts is \cite[Theorem 3.6 and Theorem 4.1]{Ri88}.

\section{Minimal Levi subgroups}\label{sec-BFM-Levi}

In the rest of this article, we will denote by $G$ a simply connected almost simple group and assume $g>1$. We explain in \S \ref{subsec-BFM-repre} how the Levi subgroups $L_z$ naturally appear in $z$-twisted character varieties. Explicit descriptions of these Levi subgroups are given in \S \ref{subsec-Table-BFM}.

\subsection{Minimal Levi subgroups and reducible representations}\label{subsec-BFM-repre}\hfill

We will fix a maximal torus $T$ as well as a Borel subgroup $B$ containing it, which determines a set of simple roots $\Delta\subset\Phi$. The Weyl group $N_G(T)/T$ will be denoted by $W$.

\begin{Lem}\label{Lem-z-in-Lz}
Let $z\in Z_G$ and let $L\subset G$ be a Levi subgroup containing $T$. Then, the following are equivalent:
\begin{itemize}
\item[(i)] The twisted representation variety $\Rep_z(\Pi,L)$ is nonempty.
\item[(ii)] $\mathbf{D}_L(z)=1\in Z_L^{\circ}/Z_L^{\circ}\cap [L,L]$.
\item[(iii)] $z\in [L,L]$.
\item[(iv)] $z\in T_1:=T\cap[L,L]$.
\item[(v)] $z$ lies in the image $\tilde{\iota}_T(\mathbb{Q}\otimes_{\mathbb{Z}}Q_L^{\vee})$.
\item[(vi)] Suppose that $L'\supset T$ is a Levi subgroup that is conjugate to $L$ and admits a basis of simple roots contained in $\Delta$. Suppose that $\varpi^{\vee}$ is the fundamental coweight corresponding to $z$ and 
$$
\varpi^{\vee}=\sum_{i=1}^{\ell}a_i\alpha_i^{\vee}.
$$
Then, $a_i\in\mathbb{Z}$ if $\alpha_i^{\vee}$ is not a coroot of $L'$.
\end{itemize}
\end{Lem}
Note that $T_1$ is a maximal torus of $[L,L]$. It is easy to see that $(T\cap[L,L])^{\circ}$ is a maximal torus of $[L,L]$. But $T\cap[L,L]$ is connected since 
$$
T\cap[L,L]\subset C_{[L,L]}((T\cap[L,L])^{\circ})=(T\cap[L,L])^{\circ}.
$$
\begin{proof}
Note that $z\in Z_L\subset T$. Recall that the map $z\mapsto\mathbf{D}_L(z)$ amounts to writing $z=z_0z_1$ with $z_0\in Z_L^{\circ}$ and $z_1\in [L,L]$ and taking the coset of $z_0$. The equivalence (ii)$\Leftrightarrow$(iii) then follows from the definition. Now we show (i)$\Leftrightarrow$(iii). If $\Rep_z(\Pi,L)$ contains a point $\rho$, then the equation defining $\Rep_z(\Pi,L)$ shows $z\in[L,L]$. Conversely, if (iii) holds, then $z\in Z_{[L,L]}$. It is known that $\Rep_z(\Pi,[L,L])$ is nonempty, and so (i) follows. Now, (iii)$\Leftrightarrow$(iv) follows from the definition of $T_1$. Next, we show (iv)$\Leftrightarrow$(v). If $z$ is a finite order element of $T_1$, then it lies in the image of $\tilde{\iota}_{T_1}$. According to \cite[Proposition 12.14]{MT}, we have $Y(T_1)\cong Q_L^{\vee}$, so that (v) holds. Since the images of the coroots of $L$ are contained in $T_1$, we have (v)$\Rightarrow$(iv). It remains to show (v)$\Leftrightarrow$(vi). Let $L'$ be a Levi subgroup as in (vi). Since $\Phi_L$ and $\Phi_{L'}$ are conjugate under $N_G(T)$ and $z$ is invariant under the conjugation by $N_G(T)$, we may replace $L$ by $L'$ in (v) and assume $L=L'$ in (vi). Suppose that (vi) holds. By construction, $z=\tilde{\iota}_T(\varpi^{\vee})$. If $a_i\in\mathbb{Z}$, then $a_i\alpha_i^{\vee}\in Q_G^{\vee}=Y(T)$, and thus is mapped to $1$ by $\tilde{\iota}_T$. We see that $z$ is the product of $\tilde{\iota}_T(a_i\alpha_i^{\vee})$ with $\alpha_i^{\vee}$ being a coroot of $L$; that is, (v) is true. Note that the kernel of $\tilde{\iota}_T$ is precisely $Y(T)$. Therefore, (v) implies that $\varpi^{\vee}$ lies in the additive group generated by $\mathbb{Q}\otimes_{\mathbb{Z}}Q_{L'}^{\vee}$ and $Y(T)=Q_G^{\vee}$, and so $a_i$ is an integer if $\alpha_i^{\vee}$ is not a coroot of $L'$.
\end{proof}

\begin{Defn}\label{Defn-Lz}
Let $z\in Z_G$ and let $\varpi^{\vee}$ be the fundamental coweight corresponding to $z$. Write 
$$
\varpi^{\vee}=\sum_{i=1}^{\ell}a_i\alpha_i^{\vee},
$$
and
$$
I(z)=\{i\mid \text{$a_i$ is not an integer}\}.
$$
Define the Levi subgroup $L_z$, which contains $T$ and has $I(z)$ as a basis of simple roots. We will write $S_z:=Z_{L_z}^{\circ}$.
\end{Defn}

\begin{Cor}
Let $I(z)\subset\Delta$ be as in Definition \ref{Defn-Lz}. If $I'\subset\Delta$ is a subset that contains $I(z)$, and $I''\subset\Delta$ is $W$-conjugate to $I'$, then $I''$ contains $I(z)$. In particular, the only subset of $\Delta$ that is $W$-conjugate to $I(z)$ is $I(z)$ itself.
\end{Cor}
\begin{proof}
The assertion follows from the fact that $z$ is fixed by $W$.
\end{proof}
\begin{Lem}\label{Lem-N_W(W_z)/W_z}
Let $L\subset G$ be a Levi subgroup containing $T$ and let $W_L$ be the Weyl group of $L$ defined by $T$. Then, we have natural isomorphisms:
$$
N_G(L)/L\cong N_W(\Phi_L)/W_L\cong N_W(W_L)/W_L.
$$
\end{Lem}
\begin{proof}
This follows from the proof of \cite[Proposition 3.4.3]{DM20}.
\end{proof}
\begin{Prop}\label{Prop-L-supset-L_z}
Let $z\in Z_G$ and let $L_z$ be the Levi subgroup given by Definition \ref{Defn-Lz}. Then, the following hold:
\begin{itemize}
\item[(i)] If $L\subset G$ is a Levi subgroup such that there exists a $z$-twisted $L$-representation, then $L$ contains a $G$-conjugate of $L_z$.
\item[(ii)] Every $z$-twisted $L_z$-representation is $L_z$-irreducible.
\end{itemize}
\end{Prop}
\begin{proof}
We first prove (i).  Up to a conjugation by $G$, we may assume that $L$ contains the given maximal torus $T$ and that a basis of $\Phi_L$ is contained in $\Delta$. It follows from the equivalence (i)$\Leftrightarrow$(vi) in Lemma \ref{Lem-z-in-Lz} and Definition \ref{Defn-Lz} that the set of the simple roots of $L$ contains that of $L_z$; hence the assertion. Now, we prove (ii). Suppose that $\rho:\Pi\rightarrow L_z$ is reducible so that $\Ima\rho\subset Q$ for some proper parabolic subgroup $Q\subset L_z$. Let $M\subset Q$ be a Levi factor of $Q$. Composing $\rho$ with the projection $Q\rightarrow M$ gives a $z$-twsited representation with values in $M$ (note that $z\in M$), which contradicts (i).
\end{proof}

\subsection{Tables for minimal Levi subgroups}\label{subsec-Table-BFM}\hfill

Let us describe the Levi subgroups $L_z$ in terms of their simple roots. Since for types $E_8$, $F_4$ and $G_2$, the simply connected $G$ is also of adjoint type, we have $Z_G=1$ and $L_z$ is a maximal torus. For other types, the set $I(z)$ of the simple roots of $L_z$ is a nonempty subset of $\Delta$. The labelling of the elements of $\Delta$ are shown in the following list of Dynkin diagrams. For each type, the highest root can be written as $\sum_{\alpha_i\in\Delta}n_i\alpha_i$, a positive linear combination of the simple roots of $G$. The coefficients $n_{\alpha}$ are indicated in the nodes of the diagrams. The elements of $Z_G$ different from $1$ are in bijection with the nodes with $n_i=1$.

$$
A_{\ell}\quad\quad\quad
\begin{tikzcd}[row sep=0em, column sep=0.2em]
\circled{1} \arrow[rr, swap, shorten=-2mm, dash, yshift=-0.1mm]  && \circled{1} \arrow[rr, swap, shorten=-2mm, dash, yshift=-0.1mm] && \cdots \arrow[rr, swap, shorten=-2mm, dash, yshift=-0.1mm] && \circled{1} \arrow[rr, swap, shorten=-2mm, dash, yshift=-0.1mm] && \circled{1}\\
\alpha_1 && \alpha_2 && \cdots && \alpha_{\ell-1} && \alpha_{\ell}
\end{tikzcd}
$$

$$
B_{\ell}\quad\quad\quad
\begin{tikzcd}[row sep=0em, column sep=0.2em]
\circled{1} \arrow[rr, swap, shorten=-2mm, dash, yshift=-0.1mm]  && \circled{2} \arrow[rr, swap, shorten=-2mm, dash, yshift=-0.1mm] && \cdots \arrow[rr, swap, shorten=-2mm, dash, yshift=-0.1mm] && \circled{2} \arrow[rr, swap, shorten=-2mm, equal, yshift=-0.1mm, ">" marking] && \circled{2}\\
\alpha_1 && \alpha_2 && \cdots && \alpha_{\ell-1} && \alpha_{\ell}
\end{tikzcd}
$$

$$
C_{\ell}\quad\quad\quad
\begin{tikzcd}[row sep=0em, column sep=0.2em]
\circled{2} \arrow[rr, swap, shorten=-2mm, dash, yshift=-0.1mm]  && \circled{2} \arrow[rr, swap, shorten=-2mm, dash, yshift=-0.1mm] && \cdots \arrow[rr, swap, shorten=-2mm, dash, yshift=-0.1mm] && \circled{2} \arrow[rr, swap, shorten=-2mm, equal, yshift=-0.1mm, "<" marking] && \circled{1}\\
\alpha_1 && \alpha_2 && \cdots && \alpha_{\ell-1} && \alpha_{\ell}
\end{tikzcd}
$$

$$
D_{\ell}\quad\quad\quad
\begin{tikzcd}[row sep=0em, column sep=0.2em]
&& && && && \circled{1} & \alpha_{\ell-1}\\
\circled{1} \arrow[rr, swap, shorten=-2mm, dash, yshift=-0.1mm]  && \circled{2} \arrow[rr, swap, shorten=-2mm, dash, yshift=-0.1mm] && \cdots \arrow[rr, swap, shorten=-2mm, dash] && \circled{2} \arrow[drr, swap, shorten=-2.5mm, dash, yshift=-0.1mm] \arrow[urr, swap, shorten=-2.5mm, dash, yshift=-0.1mm] &&  &\\
\alpha_1 && \alpha_2 && \cdots && \alpha_{\ell-2} && \circled{1} & \alpha_{\ell} 
\end{tikzcd}
$$

$$
E_6\quad\quad\quad
\begin{tikzcd}[row sep=0em, column sep=-0.1em]
&{}&& \!\!\alpha_2\!\!\!\! & \circled{2} \arrow[dd, swap, shorten=-1.5mm, dash] &{}& &{}&\\
&&&&{}&&&&\\
\circled{1} \arrow[rr, swap, shorten=-2mm, dash, yshift=-0.1mm]  && \circled{2} \arrow[rr, swap, shorten=-2mm, dash, yshift=-0.1mm] && \circled{3} \arrow[rr, swap, shorten=-2mm, dash, yshift=-0.1mm] && \circled{2} \arrow[rr, swap, shorten=-2mm, dash, yshift=-0.1mm] && \circled{1}\\
\alpha_1 && \alpha_3 && \alpha_4 && \alpha_{5} && \alpha_6
\end{tikzcd}
$$

$$
E_7\quad\quad\quad
\begin{tikzcd}[row sep=0em, column sep=-0.1em]
&{}&& \!\!\alpha_2\!\!\!\! & \circled{2} \arrow[dd, swap, shorten=-1.5mm, dash] &{}& &{}& &{}&\\
&&&&{}&&&&&&\\
\circled{2} \arrow[rr, swap, shorten=-2mm, dash, yshift=-0.1mm]  && \circled{3} \arrow[rr, swap, shorten=-2mm, dash, yshift=-0.1mm] && \circled{4} \arrow[rr, swap, shorten=-2mm, dash, yshift=-0.1mm] && \circled{3} \arrow[rr, swap, shorten=-2mm, dash, yshift=-0.1mm] && \circled{2} \arrow[rr, swap, shorten=-2mm, dash, yshift=-0.1mm] && \circled{1} \\
\alpha_1 && \alpha_3 && \alpha_4 && \alpha_{5} && \alpha_6 && \alpha_7
\end{tikzcd}
$$

For any $z\ne 1$, let $\varpi^{\vee}$ be the fundamental coweight corresponding to $z$. Write $\varpi^{\vee}=\sum_{i=1}^{\ell}a_i\alpha_i^{\vee}$. The coefficients $a_i$ are given by the inverse of the Cartan matrix, which can be found in the book of Onishchik-Vinberg \cite[Reference Chapter, \S 2, Table 2]{OV}. If $\omega^{\vee}=\omega^{\vee}_j$ for some $j\in\{1,2,\ldots,\ell\}$, then the coefficients $\{a_1,\ldots,a_{\ell}\}$ are given by the $j$-th row of that table. Note that for exceptional types, the labelling of $\Delta$ in \cite{OV} is different from the labelling above. 

\begin{center}
\begin{tabular}{|c|c|c|c|c|}
\hline
Type & $Z_G$ & $z\in Z_G$ & $I(z)$ & $|I(z)|$ \\
\hline
$A_{\ell}$ & $\mathbb{Z}/(\ell+1))$ & $m$ & $\Delta\setminus\{\alpha_i\mid\text{ $(\ell+1)$ divides $m(\ell+1-i)$}\}$ & $\gcd(\ell+1,m)-1$ \\
\hline
$B_{\ell}$ & $\mathbb{Z}/2$ & $1$ & $\{\alpha_{\ell}\}$ & 1\\
\hline
$C_{2m}$ & $\mathbb{Z}/2$ & $1$ & $\{\alpha_1,\alpha_3,\ldots,\alpha_{\ell-1}\}$ & $m$ \\
\hline
$C_{2m+1}$ & $\mathbb{Z}/2$ & $1$ & $\{\alpha_1,\alpha_3,\ldots,\alpha_{\ell} \}$ & $m+1$\\
\hline
$D_{\ell}$ & $\mathbb{Z}/4$ or $(\mathbb{Z}/2)^2$ & $\alpha_1$ & $\{\alpha_{\ell-1},\alpha_{\ell}\}$ & 2 \\
\hline
$D_{2m+1}$ & $\mathbb{Z}/4$ & $\alpha_{\ell-1}$ or $\alpha_{\ell}$ & $\{\alpha_1,\alpha_3,\ldots,\alpha_{\ell-2},\alpha_{\ell-1},\alpha_{\ell}\}$ & $m+2$\\
\hline
$D_{4m}$ & $(\mathbb{Z}/2)^2$ & $\alpha_{\ell-1}$ & $\{\alpha_1,\alpha_3,\ldots,\alpha_{\ell-3},\alpha_{\ell}\}$ & $2m$ \\
\hline
$D_{4m}$ & $(\mathbb{Z}/2)^2$ & $\alpha_{\ell}$ & $\{\alpha_1,\alpha_3,\ldots,\alpha_{\ell-3},\alpha_{\ell-1}\}$ & $2m$ \\
\hline
$D_{4m+2}$ & $(\mathbb{Z}/2)^2$ & $\alpha_{\ell-1}$ & $\{\alpha_1,\alpha_3,\ldots,\alpha_{\ell-3},\alpha_{\ell-1}\}$ & $2m+1$ \\
\hline
$D_{4m+2}$ & $(\mathbb{Z}/2)^2$ & $\alpha_{\ell}$ & $\{\alpha_1,\alpha_3,\ldots,\alpha_{\ell-3},\alpha_{\ell}\}$ & $2m+1$ \\
\hline
$E_6$ & $\mathbb{Z}/3$ & $1$ or $2$ & $\{\alpha_1,\alpha_3,\alpha_5,\alpha_6\}$ & 4\\
\hline
$E_7$ & $\mathbb{Z}/2$ & $1$ & $\{\alpha_2,\alpha_5,\alpha_7\}$ & 3\\
\hline
\end{tabular}
\end{center}
In the above table, the element $z$ is represented either by a vertex of the Dynkin diagram or by a number in an (additive) cyclic group; however, an element of $Z_G$ is often written multiplicatively in the rest of this article.

\section{Stabilisers}

This section begins with the definition of the subgroups $H_z\subset L_z$ as in Theorem \ref{Thm-B} and the computation of their connected centralisers $G_z^{\circ}$. For a $z$-twisted $H_z$-valued representation $\rho$, its connected stabiliser $\Stab_G(\rho)^{\circ}$ may be strictly larger than $G_z^{\circ}$ if $G$ is of type $B$ or $C$. We will see in \S \ref{subsec-type-BC} that this only happens for $G$ of type $C_{2m+1}$, where we define a subgroup $H_z'\subset H_z$ with a larger connected centraliser $G_z^{\prime\circ}$. We then determine the connected components of $G_z$ and $G_z'$ in \S \ref{subsec-Gamma-G_z}. In \S \ref{subsec-large-stab}, we will show that $G_z$ and $G_z'$ are indeed the largest possible stabiliser groups.

\subsection{$H_z$ and $C_G(H_z)^{\circ}$}\label{subsec-Defn-H_z}\hfill

We define in each case a particular irreducible subgroup $H_z\subset L_z$ containing $z$. 

\textit{Type $A_{\ell}$}. Suppose that $z\in Z_G$ is of order $d$. The derived subgroup $[L_z,L_z]$ is isomorphic to the direct product of $(\ell+1)/d$ copies of $\SL_d$. We take $H_z=\SL_d$ diagonally embedded in $[L_z,L_z]$.

\textit{Type $B_{\ell}$}. The derived subgroup $[L_z,L_z]$ is isomorphic to $\SL_2$. We take $H_z=[L_z,L_z]$.

\textit{Type $C_{2m}$.} The derived subgroup $[L_z,L_z]$ is isomorphic to a direct product of $m$-copies of $\SL_2$. We take for $H_z$ the $\SL_2$ diagonally embedded in $[L_z,L_z]$. 

\textit{Type $C_{2m+1}$.} The derived subgroup $[L_z,L_z]$ is isomorphic to a direct product of $(m+1)$-copies of $\SL_2$, the last factor having the long simple root as its positive root. We take $H_z\cong\SL_2\times\SL_2$, the first factor being diagonally embedded in the first $m$ $\SL_2$-factors of $[L_z,L_z]$ and the second being a copy of the $\SL_2$-factor corresponding to the long root.  

\textit{Type $D_{\ell}$.} Suppose that $z\in Z_G$ corresponds to $\alpha_1$. The derived subgroup $[L_z,L_z]\cong\SL_2\times\SL_2$ and the two positive roots are $\alpha_{\ell-1}$ and $\alpha_{\ell}$. We take for $H_z$ the $\SL_2$ diagonally embedded in $[L_z,L_z]$.

\textit{Type $D_{2m+1}$.} Suppose that $z\in Z_G$ corresponds to $\alpha_{\ell-1}$ or $\alpha_{\ell}$. The derived subgroup $[L_z,L_z]$ is isomorphic to the direct product of $m-1$ copies of $\SL_2$ and one $\SL_4$; the $\SL_4$-factor has a basis of simple roots consisting of $\alpha_{\ell-2}$, $\alpha_{\ell-1}$ and $\alpha_{\ell}$. We take $H_z=\SL_2\times\SL_4$ with the first factor diagonally embedded in the $m-1$ copies of $\SL_2$ in $[L_z,L_z]$.

\textit{Type $D_{4m}$.} Suppose that $z\in Z_G$ corresponds to $\alpha_{\ell}$. The derived subgroup $[L_z,L_z]$ is isomorphic to a direct product of $2m$ copies of $\SL_2$. We take for $H_z$ the $\SL_2$ diagonally embedded in $[L_z,L_z]$. Other type $D$ groups of even ranks are similar.

\textit{Type $E_6$.} The derived subgroup $[L_z,L_z]$ is isomorphic to $\SL_3\times\SL_3$, which has $\{\alpha_1,\alpha_3,\alpha_5,\alpha_6\}$ as a basis of simple roots. We take $H_z=\SL_3$ diagonally embedded in $[L_z,L_z]$.

\textit{Type $E_7$.} The derived subgroup $[L_z,L_z]$ is isomorphic to $\SL_2\times\SL_2\times\SL_2$, which has positive roots $\{\alpha_2,\alpha_5,\alpha_7\}$. We take $H_z=\SL_2$ diagonally embedded in $[L_z,L_z]$.

Note that $H_z$ is an irreducible subgroup of $[L_z,L_z]$, since there is no Levi subgroup of the latter containing $H_z$.

In each case, we write $[L_z,L_z]\cong \prod_{i\in I}G_i$, where every $G_i$ is a special linear group. The subgroup $H_z\subset[L_z,L_z]$ is the fixed point subgroup of automorphism $\sigma$ cyclically permuting some direct factors of $[L_z,L_z]$, say $\{G_i\mid i\in \tilde{I}\}$ for some subset $\tilde{I}\subset I$. 
\begin{Lem}\label{Lem-L'-SL}
The subgroup $\prod_{i\in \tilde{I}}G_i\subset[L_z,L_z]$ is always contained in a Levi subgroup $L'\subset G$ whose derived subgroup is a special linear group. In particular, the cyclic permutation $\sigma$ defining $H_z$ can be realised by some $n\in[L',L']$. 
\end{Lem}
\begin{proof}
It suffices to notice that in each case, there exists a subset $I'\subset \Delta$ containing $\tilde{I}$, which is the vertices of a connected Dynkin subdiagram of type $A$.

\textit{Type $A_{\ell}$}. We take $I'=\Delta$ in this case.

\textit{Type $B_{\ell}$}. In this case, the derived subgroup $[L_z,L_z]$ is equal to $H_z$ and the cyclic permutation defining $H_z$ is trivial.

\textit{Type $C_{2m}$.} Let $I'$ be the set of short simple roots in $\Delta$, which is the base of a root system of type $A_{2m-1}$. The set $I(z)=\tilde{I}$ is a subset of $I'$. 

\textit{Type $C_{2m+1}$.} Let $I'=\Delta\setminus\{\alpha_{\ell},\alpha_{\ell-1}\}$, which is the base of a root system of type $A_{2m-1}$. The set of short simple roots in $I(z)$ is precisely $\tilde{I}$. We have $\tilde{I}\subset I'$. 

\textit{Type $D_{\ell}$.} Suppose that $z\in Z_G$ corresponds to $\alpha_1$. We have $\tilde{I}=I(z)=\{\alpha_{\ell-1},\alpha_{\ell}\}$. We take $I'=\{\alpha_{\ell-2},\alpha_{\ell-1},\alpha_{\ell}\}$, which spans a root system of type $A_3$.

\textit{Type $D_{2m+1}$.} Suppose that $z\in Z_G$ corresponds to $\alpha_{\ell-1}$ or $\alpha_{\ell}$. In this case, the set $I(z)$ is the union of $\tilde{I}$ and $\{\alpha_{\ell-2},\alpha_{\ell-1},\alpha_{\ell}\}$, the latter of which forms a Dynkin diagram of type $A_3$. We take 
$$
I'=\Delta\setminus\{\alpha_{\ell-3},\alpha_{\ell-2},\alpha_{\ell-1},\alpha_{\ell}\}.
$$

\textit{Type $D_{4m}$.} Suppose that $z\in Z_G$ corresponds to $\alpha_{\ell}$. In this case, we have 
$$
\tilde{I}=I(z)=\{\alpha_1,\alpha_3,\ldots,\alpha_{\ell-3},\alpha_{\ell-1}\}.
$$
We take $I'=\Delta\setminus\{\alpha_{\ell}\}$, which is of type $A_{\ell-1}$.

\textit{Type $E_6$.} We have $\tilde{I}=I(z)=\{\alpha_1,\alpha_3,\alpha_5,\alpha_6\}$ as a basis of simple roots. We take $I'=\Delta\setminus\{\alpha_2\}$, which is of type $A_5$.

\textit{Type $E_7$.} We have $\tilde{I}=I(z)=\{\alpha_2,\alpha_5,\alpha_7\}$. We take $I'=\Delta\setminus\{\alpha_1,\alpha_3\}$, which is of type $A_5$.
\end{proof}
\begin{Rem}
It follows from the above lemma that $z\in H_z$. By definition, we have $H_z=C_{[L_z,L_z]}(n)$. However, $z$ is a central element of $G$, and thus it commutes with $n$.
\end{Rem}

Let $H\subset L_z$ be an irreducible subgroup. Note that $H_z$ is such a group. Lemma \ref{Lem-S<C_G(H)} shows that $C_G(H)^{\circ}$ contains $S_z$ as a maximal torus.  
\begin{Lem}\label{Lem-W{circ}->N(L)/L}
There is a natural embedding
$$
W_{C_G(H)^{\circ}}:=N_{C_G(H)^{\circ}}(S_z)/S_z\hooklongrightarrow N_G(L_z)/L_z.
$$
\end{Lem}
\begin{proof}
Let $n\in N_{C_G(H)^{\circ}}(S_z)$. Since $L_z=C_G(S_z)$, we have $n\in N_G(L_z)$. This shows that we have an inclusion $N_{C_G(H)^{\circ}}(S_z)\hookrightarrow N_G(L_z)$. The elements of $L_z$ acts trivially on $S_z$. An element of $N_{C_G(H)^{\circ}}(S_z)$ acts trivially on $S_z$ if and only if it lies in $S_z$. This shows the injectivity of the desired embedding.
\end{proof}
\begin{Lem}\label{Lem-W->N(L)/L}
There is a natural homomorphism
$$
W_{C_G(H)}:=N_{C_G(H)}(S_z)/S_z\longrightarrow N_G(L_z)/L_z,
$$
and its kernel is isomorphic to $C_{L_z}(H)/S_z\cong Z_{L_z}/S_z$. In particular, if $Z_{L_z}=S_z$ (i.e., $Z_{L_z}$ is connected), then this homomorphism is injective.
\end{Lem}
\begin{proof}
The well-definedness of this homomorphism follows from the same arguments for the above lemma. Since $L_z$ automatically normalises $S_z$, the kernel is $(L_z\cap C_G(H))/S_z$. The inclusion $Z_{L_z}\subset C_{L_z}(H)$ is immediate. Let $\tilde{H}\subset S_z\times [L_z,L_z]$ be the inverse image of $H$, and let $H_1\subset[L_z,L_z]$ be its projection to $[L_z,L_z]$. Then, $H_1\subset[L_z,L_z]$ is an irreducible subgroup. Note that $[L_z,L_z]$ is a direct product of special linear groups. The projection of $H_1$ to each direct factor of $[L_z,L_z]$ is necessarily an irreducible subgroup. For any $r>1$, if $H'\subset \SL_r$ is an irreducible subgroup, then $C_{\SL_r}(H')=Z_{\SL_r}$, which is simply Schur's lemma. Therefore, we have $C_{[L_z,L_z]}(H)=C_{[L_z,L_z]}(H_1)=Z_{[L_z,L_z]}$, which is contained in $Z_{L_z}$, and thus $C_{L_z}(H)\subset Z_{L_z}$, as desired.
\end{proof}
\begin{Lem}
Let $\mathfrak{s}_z$ be the Lie algebra of $S_z$, which is a subspace of $\mathfrak{t}$. Then, there is a natural isomorphism
$$
N_G(L_z)/L_z\lisom N_W(\mathfrak{s}_z)/C_W(\mathfrak{s}_z),
$$
where $C_W(\mathfrak{s}_z)$ is the subgroup of $W$ consisting of the elements that fix $\mathfrak{s}_z$ pointwise and $N_W(\mathfrak{s}_z)$ consists of the elements that preserve $\mathfrak{s}_z$.
\end{Lem}
\begin{proof}
By Lemma \ref{Lem-N_W(W_z)/W_z}, it suffices to give an isomorphism
$$
N_W(W_L)/W_L\lisom N_W(\mathfrak{s}_z)/C_W(\mathfrak{s}_z).
$$
We have $C_W(\mathfrak{s}_z)=W_L$, since both groups are the subgroup generated by the reflections that fix $\mathfrak{s}_z$ pointwise. It follows that if $w\in W$ normalises $\mathfrak{s}_z$, then it normalises $W_L$. The inclusion $N_W(W_L)\subset N_W(\mathfrak{s}_z)$ follows from the observation that $\mathfrak{s}_z\subset\mathfrak{t}$ is precisely the fixed point subspace of $W_L$. 
\end{proof}
\begin{Rem}\label{Rem-type-Phi-proj}
The group $N_W(\mathfrak{s}_z)/C_W(\mathfrak{s}_z)$ can be identified with the Weyl group of certain root system defined on $\mathfrak{s}_z$. Note that $\mathfrak{s}_z=\mathfrak{t}^{w_z}$, where $w_z\in W$ is the element preserving $\Delta$ and corresponds to $z$ (see \cite[Equation (3.7)]{Bon}); this is a combination of \cite[Proposition 3.5.4]{BFM} and \cite[Proposition 4.2.1]{BFM}. Applying \cite[Proposition 6.3.4 and Proposition 6.5.5]{BFM} to $V:=\mathfrak{t}$ and $l:=w_z$, we deduce that $N_W(\mathfrak{s}_z)/C_W(\mathfrak{s}_z)$ is isomorphic to the Weyl group of the root system $\Phi^{proj}(l)$. For every type of $G$ and every nontrivial $z\in Z_G$ (the group generated by $z$ is denoted by $\mathcal{C}$ in \cite{BFM}), the corresponding type of $\Phi^{proj}$ is given in the table Root Systems on $\mathfrak{t}^{w_{\mathcal{C}}}$ in \cite[Diagrams and tables]{BFM}. For the convenience of the reader, we list the type of the Weyl group of $\Phi^{proj}$ as follows.
\begin{center}
\begin{tabular}{|c|c|c|c|c|c|c|c|c|c|}
\hline
$G$ & $A_{\ell}$ & $B_{\ell}$ & $C_{2m}$ & $C_{2m+1}$ & $D_{\ell}$ & $D_{2m+1}$ & $D_{2m}$ & $E_6$ & $E_7$ \\
\hline
$z$ & order $d$ & $\ne 1$ & $\ne 1$ & $\ne1$ & $\alpha_1$ & $\alpha_{\ell-1}$ or $\alpha_{\ell}$ & $\alpha_{\ell-1}$ or $\alpha_{\ell}$ & $\ne1$ & $\ne1$ \\ 
\hline
$W(\Phi^{proj})$ & $A_{d-1}$ & $B_{\ell-1}$ & $B_{m}$ & $B_{m}$ & $B_{\ell-2}$ & $B_{m-1}$ & $B_{m-1}$ & $G_2$ & $F_4$ \\
\hline
\end{tabular}
\end{center}
\end{Rem}

\begin{Thm}\label{Thm-Cent-H_z}
The connected centralisers $C_G(H_z)^{\circ}$ are as follows.
\begin{itemize}
\item Type $A_{\ell}$; $z$ is of order $d\mid \ell+1$.
$$
C_G(H_z)^{\circ}\cong\SL_{(\ell+1)/d}.
$$
\item Type $B_{\ell}$. 
$$
C_G(H_z)^{\circ}\cong\Spin_{2\ell-2}.
$$
\item Type $C_{2m}$ and Type $C_{2m+1}$. 
$$
C_G(H_z)^{\circ}\cong\SO_{2m}.
$$
\item Type $D_{\ell}$; $z\in Z_G$ corresponds to $\alpha_1$. 
$$
C_G(H_z)^{\circ}\cong\Spin_{2\ell-3}.
$$
\item Type $D_{2m+1}$; $z\in Z_G$ corresponds to $\alpha_{\ell-1}$ or $\alpha_{\ell}$. 
$$
C_G(H_z)^{\circ}\cong\Sp_{2(m-1)}.
$$
\item Type $D_{4m}$ and type $D_{4m+2}$; $z\in Z_G$ corresponds to $\alpha_{\ell-1}$ or $\alpha_{\ell}$. 
$$
C_G(H_z)^{\circ}\cong\Sp_{4m}.
$$
\item Type $E_6$. 
$$
C_G(H_z)^{\circ} \textit{ is semi-simple of type $G_2$.}
$$
\item Type $E_7$. 
$$
C_G(H_z)^{\circ} \textit{ is semi-simple of type $F_4$.}
$$
\end{itemize}
\end{Thm}
\begin{proof}
The root system of $C_G(H_z)^{\circ}$ will be determined in the following manner. The Lie algebra $\mathfrak{g}:=\Lie G$ decomposes into the weight spaces of $S_z$, which is a maximal torus of $C_G(H_z)^{\circ}$. For any $\chi\in X(S_z)$, we consider $\mathfrak{g}(\chi)$, the subspace of $\mathfrak{g}$ on which $S_z$ acts via $\chi$, which is a representation of $H_z$, since $H_z$ commutes with $S_z$. Note that the Lie algebra of $C_G(H_z)^{\circ}$ is the space $\mathfrak{g}^{H_z}$ of invariant vectors. If the subspace $\mathfrak{g}(\chi)^{H_z}$ is nonzero, then $\chi$ is a root of $C_G(H_z)^{\circ}$.

\textit{Type $A_{\ell}$.} This case can be treated without any Lie-theoretic computations. The asserted result follows from the fact that the automorphism group of the direct product of $(\ell+1)/d$ copies of the standard representation of $\SL_d$ is isomorphic to $\GL_{(\ell+1)/d}$.

\textit{Type $B_{\ell}$.} Let $\{\epsilon_i\}_{1\le i\le \ell}$ be a basis of $X(T)$ so that the simple roots are $\alpha_i=\epsilon_i-\epsilon_{i+1}$ for $1\le i\le \ell-1$ and $\alpha_{\ell}=\epsilon_{\ell}$. Note that the positive root of $H_z\cong\SL_2$ is $\alpha_{\ell}$. For $1\le i<j\le \ell-1$, the root subspace $\mathfrak{g}_{\pm\epsilon_i\pm\epsilon_j}$ is a one-dimensional representation of $H_z$, and thus it is contained in $\mathfrak{g}^{H_z}$. The subspace spanned by $\{\mathfrak{g}_{\epsilon_i-\epsilon_{\ell}},\mathfrak{g}_{\epsilon_{i}},\mathfrak{g}_{\epsilon_i+\epsilon_{\ell}}\}$ forms an irreducible representation of $H_z$ and there is no nonzero invariant vector in this representation. This shows that $C_G(H_z)^{\circ}$ has a root system of type $D_{\ell-1}$. The cocharacter lattice $Y(S_z)$ consists of the vectors in $Y(T)=Q_G^{\vee}$ on which $\epsilon_{\ell}$ vanishes, and we find that it is spanned by the coroots $\pm\epsilon_i\pm\epsilon_j$ for $1\le i<j\le \ell-1$. It follows that $C_G(H_z)^{\circ}\cong\Spin_{2\ell-2}$. 

\textit{Type $C_{2m}$.} Let $\{\epsilon_i\}_{1\le i\le \ell}$ be a basis of $X(T)$ so that the simple roots are $\alpha_i=\epsilon_i-\epsilon_{i+1}$ for $1\le i\le \ell-1$ and $\alpha_{\ell}=2\epsilon_{\ell}$. The positive roots of $L_z$ are $\{\epsilon_{2j-1}-\epsilon_{2j}\}_{1\le j\le m}$. For any $1\le j<k\le m$, the span of 
$$
\{\mathfrak{g}_{\epsilon_{2j}-\epsilon_{2k-1}},\mathfrak{g}_{\epsilon_{2j-1}-\epsilon_{2k-1}},\mathfrak{g}_{\epsilon_{2j}-\epsilon_{2k}},\mathfrak{g}_{\epsilon_{2j-1}-\epsilon_{2k}}\}
$$
forms a representation of $L_z$, and similarly for
$$
\{\mathfrak{g}_{\epsilon_{2j}+\epsilon_{2k}},\mathfrak{g}_{\epsilon_{2j-1}+\epsilon_{2k}},\mathfrak{g}_{\epsilon_{2j}+\epsilon_{2k-1}},\mathfrak{g}_{\epsilon_{2j-1}+\epsilon_{2k-1}}\}.
$$
The coroot of $H_z$ is $\sum_{1\le j\le m}(\epsilon_{2j-1}-\epsilon_{2j})$, with which we can determine the weights of the action of $H_z$ on the above representations. It turns out that each of the above representation decomposes into the direct sum of a one-dimensional representation and the adjoint representation of $\SL_2$ and thus contains a root subspace of $C_G(H_z)^{\circ}$. The representation spanned by 
$$
\{\mathfrak{g}_{2\epsilon_{2j-1}},\mathfrak{g}_{\epsilon_{2j-1}+\epsilon_{2j}},\mathfrak{g}_{2\epsilon_{2j}}\}
$$
does not contain any $H_z$-invariant vector. We see that the root system of $C_G(H_z)^{\circ}$ is of type $D_{m}$. The cocharacter lattice $Y(S_z)$ is spanned by $\{e_j:=\epsilon_{2j-1}+\epsilon_{2j}\}_{1\le j\le m}$, and the coroots of $C_G(H_z)^{\circ}$ are $\pm e_j\pm e_k$ for $1\le j<k\le m$. This means that $C_G(H_z)^{\circ}\cong\SO_{2m}$.

\textit{Type $C_{2m+1}$.} The positive roots of $L_z$ are $\{\epsilon_{2j-1}-\epsilon_{2j}\}_{1\le j\le m}$ and $2\epsilon_{\ell}$. The same construction above produces a root system of type $D_m$ since the second $\SL_2$-factor of $H_z$ acts trivially on the relevant weight subspaces of form $\mathfrak{g}(\chi)$. There are additional representations spanned by $$
\{\mathfrak{g}_{\epsilon_{2j}-\epsilon_{\ell}},\mathfrak{g}_{\epsilon_{2j-1}-\epsilon_{\ell}},\mathfrak{g}_{\epsilon_{2j}+\epsilon_{\ell}},\mathfrak{g}_{\epsilon_{2j-1}+\epsilon_{\ell}}\}
$$
which contains no nonzero invariant vector of $H_z$. Again, we have $C_G(H_z)^{\circ}\cong\SO_{2m}$.

\textit{Type $D_{\ell}$.} Suppose that $z\in Z_G$ corresponds to $\alpha_1$. For $1\le j<k\le \ell-2$, the root subspaces $\mathfrak{g}_{\pm\epsilon_j\pm\epsilon_k}$ are trivial representations of $H_z\cong\SL_2$, and for $1\le j\le \ell-2$ the span of the root subspaces
$$
\{\mathfrak{g}_{\epsilon_j-\epsilon_{\ell-1}},\mathfrak{g}_{\epsilon_j-\epsilon_{\ell}},\mathfrak{g}_{\epsilon_j+\epsilon_{\ell-1}},\mathfrak{g}_{\epsilon_j+\epsilon_{\ell}}\}
$$
forms a representation of $L_z$. The coroot of $H_z$ is $\alpha_{\ell-1}^{\vee}+\alpha_{\ell}^{\vee}=2\epsilon_{\ell-1}$. The weights of the above representation show that it contains a one dimensional trivial representation. We see that $C_G(H_z)^{\circ}$ is of type $B_{\ell-2}$, since the roots $\epsilon_j\pm\epsilon_{\ell-1}$ and $\epsilon_j\pm\epsilon_{\ell}$ restrct to a short root on $S_z$. The lattice $Y(S_z)$ is spanned by $\{\pm\epsilon_j\pm\epsilon_k\}$ for $1\le j<k\le \ell-2$ and thus coincides with the coroot lattice of type $B_{\ell-2}$. Therefore, we have $C_G(H_z)^{\circ}\cong \Spin_{2\ell-3}$. 

\textit{Type $D_{2m+1}$.} Suppose that $z\in Z_G$ corresponds to $\alpha_{\ell-1}$ or $\alpha_{\ell}$. For $1\le i\le 2m-2$, the root subspaces
\begin{equation}\label{eq-Thm-Cent-H_z-D2m+1}
\{\mathfrak{g}_{\epsilon_i-\epsilon_{\ell-2}},\mathfrak{g}_{\epsilon_i-\epsilon_{\ell-1}},\mathfrak{g}_{\epsilon_i-\epsilon_{\ell}},\mathfrak{g}_{\epsilon_i+\epsilon_{\ell}},\mathfrak{g}_{\epsilon_i+\epsilon_{\ell-1}},\mathfrak{g}_{\epsilon_i+\epsilon_{\ell-2}}\}
\end{equation}
span a six dimensional irreducible representation of the $\SL_4$-factor of $[L_z,L_z]$, and the subspace $\mathfrak{g}_{\epsilon_i+\epsilon_{\ell-2}}$ consists of highest weight vectors. Therefore, there is no nonzero $H_z$-invariant vectors in this space. For $1\le j<k\le m-1$, the root subspaces
$$
\{\mathfrak{g}_{\epsilon_{2j}-\epsilon_{2k-1}},\mathfrak{g}_{\epsilon_{2j-1}-\epsilon_{2k-1}},\mathfrak{g}_{\epsilon_{2j}-\epsilon_{2k}},\mathfrak{g}_{\epsilon_{2j-1}-\epsilon_{2k}}\}
$$
span a representation of $L_z$, and similarly for
$$
\{\mathfrak{g}_{\epsilon_{2j}+\epsilon_{2k}},\mathfrak{g}_{\epsilon_{2j-1}+\epsilon_{2k}},\mathfrak{g}_{\epsilon_{2j}+\epsilon_{2k-1}},\mathfrak{g}_{\epsilon_{2j-1}+\epsilon_{2k-1}}\}.
$$
The coroot of the $\SL_2$-factor of $H_z$ is $\sum_{1\le j\le m-1}(\epsilon_{2j-1}-\epsilon_{2j})$. Weight computation shows that such a representation contains a one dimensional trivial representation. The $\SL_4$-factor acts trivially on these spaces. The subspaces $\mathfrak{g}_{\epsilon_{2j-1}+\epsilon_{2j}}$ are also trivial representations of $H_z$. The cocharacter lattice $Y(S_z)$ is spanned by $\{\epsilon_{2j-1}+\epsilon_{2j}\}$ for $1\le j\le m-1$. Therefore, the weights of the $S_z$-action on these trivial representations form a root system of type $C_{m-1}$.

\textit{Type $D_{4m}$.} Suppose that $z\in Z_G$ corresponds to $\alpha_{\ell}$. For $1\le j<k\le 2m$, the root subspaces
$$
\{\mathfrak{g}_{\epsilon_{2j}-\epsilon_{2k-1}},\mathfrak{g}_{\epsilon_{2j-1}-\epsilon_{2k-1}},\mathfrak{g}_{\epsilon_{2j}-\epsilon_{2k}},\mathfrak{g}_{\epsilon_{2j-1}-\epsilon_{2k}}\}
$$
span a representation of $L_z$, and similarly for
$$
\{\mathfrak{g}_{\epsilon_{2j}+\epsilon_{2k}},\mathfrak{g}_{\epsilon_{2j-1}+\epsilon_{2k}},\mathfrak{g}_{\epsilon_{2j}+\epsilon_{2k-1}},\mathfrak{g}_{\epsilon_{2j-1}+\epsilon_{2k-1}}\}.
$$
Each of them contains a one dimensional trivial representation. The subspaces $\mathfrak{g}_{\epsilon_{2j-1}+\epsilon_{2j}}$ are also trivial representations of $L_z$. This gives a root system of type $C_{2m}$. Computing $Y(S_z)$ shows that it is spanned by the coroots of this root system, which implies $C_G(H_z)^{\circ}\cong\Sp_{4m}$. Similarly, if $z\in Z_G$ corresponds to $\alpha_{\ell-1}$, we again have $C_G(H_z)^{\circ}\cong\Sp_{4m}$. 

\textit{Type $E_6$.} Suppose that $z\in Z_G$ is not the identity. Let $\{\epsilon_i\}_{1\le i\le 8}$ be a standard orthonormal basis of $\mathbb{Q}^8$ so that $Y(T)\otimes\mathbb{Q}$ is the subspace orthogonal to $\epsilon_6-\epsilon_7$ and $\epsilon_7+\epsilon_8$. The positive roots are $\pm\epsilon_i+\epsilon_j$ for $1\le i<j\le 5$ and $\frac{1}{2}(\epsilon_8-\epsilon_7-\epsilon_6+\sum_{i=1}^5\nu_i\epsilon_i)$ with $\nu_i\in\{\pm1\}$ satisfying $\prod_{i=1}^5\nu_i=1$, and the simple roots are 
$$
\alpha_1=\frac{1}{2}(\epsilon_8-\epsilon_7-\epsilon_6+\epsilon_1-\epsilon_2-\epsilon_3-\epsilon_4-\epsilon_5)
$$
and 
$$
\alpha_2=\epsilon_1+\epsilon_2,\quad \alpha_3=\epsilon_2-\epsilon_1,\quad\alpha_4=\epsilon_3-\epsilon_2,\quad\alpha_5=\epsilon_4-\epsilon_3,\quad\alpha_6=\epsilon_5-\epsilon_4.
$$
The simple roots of $L_z$ are $\{\alpha_1,\alpha_3,\alpha_5,\alpha_6\}$. The root subspace $\mathfrak{g}_{\alpha_2}$ forms a one-dimensional representation of $L_z$, and thus it gives rise to a root of $C_G(H_z)$. The root $\alpha_4$ and the simple roots of $L_z$ form a root system of type $A_5$; therefore, there is a nine-dimensional nilpotent Lie subalgebra of $\mathfrak{g}$ on which $L_z$ acts with highest weight $\alpha_1+\alpha_3+\alpha_4+\alpha_5+\alpha_6$ and lowest weight $\alpha_4$. This representation is the action of $\SL_3\cong H_z$ on $\gl_3$ by conjugation, and thus contains a one-dimensional subspace of invariant vectors. The restriction of $\alpha_2$ (resp. $\alpha_4$) to $Y(S_z)\otimes\mathbb{Q}$ is the short (resp. long) simple root of a root system of type $G_2$. Moreover, it can be checked that the nonzero weights for the action of $S_z$ on $\mathfrak{g}$ are precisely the roots of this root system. It follows that $C_G(H_z)^{\circ}$ is the semi-simple group of type $G_2$.

\textit{Type $E_7$.} Suppose that $z\in Z_G$ is not the identity. Let $\{\epsilon_i\}_{1\le i\le 8}$ be a standard orthonormal basis of $\mathbb{Q}^8$ so that $Y(T)\otimes\mathbb{Q}$ is the subspace orthogonal to $\epsilon_7+\epsilon_8$. The positive roots are $\epsilon_8-\epsilon_7$, $\pm\epsilon_i+\epsilon_j$ for $1\le i<j\le 6$ and $\frac{1}{2}(\epsilon_8-\epsilon_7+\sum_{i=1}^6\nu_i\epsilon_i)$ with $\nu_i\in\{\pm1\}$ satisfying $\prod_{i=1}^6\nu_i=-1$, and the simple roots are 
$$
\alpha_1=\frac{1}{2}(\epsilon_8-\epsilon_7+\epsilon_1-\epsilon_2-\epsilon_3-\epsilon_4-\epsilon_5-\epsilon_6)
$$
and 
$$
\alpha_2=\epsilon_1+\epsilon_2,\quad \alpha_3=\epsilon_2-\epsilon_1,\quad\alpha_4=\epsilon_3-\epsilon_2,\quad\alpha_5=\epsilon_4-\epsilon_3,\quad\alpha_6=\epsilon_5-\epsilon_4,\quad\alpha_7=\epsilon_6-\epsilon_5.
$$
The simple roots of $L_z$ are $\{\alpha_2,\alpha_5,\alpha_7\}$. The root subspaces $\mathfrak{g}_{\alpha_1}$ and $\mathfrak{g}_{\alpha_3}$ are one-dimensional representations of $L_z$, and thus they give rise to roots of $C_G(H_z)$. The root subspaces 
$$
\{\mathfrak{g}_{\alpha_4},\mathfrak{g}_{\alpha_2+\alpha_4},\mathfrak{g}_{\alpha_4+\alpha_5},\mathfrak{g}_{\alpha_2+\alpha_4+\alpha_5}\}
$$ 
span a four-dimensional representation of $H_z\cong \SL_2$, which is isomorphic to the adjoint representation of $\SL_2$ on $\gl_2$; therefore, there is a one-dimensional subspace of invariant vectors. The same is true for 
$$
\{\mathfrak{g}_{\alpha_6},\mathfrak{g}_{\alpha_6+\alpha_5},\mathfrak{g}_{\alpha_6+\alpha_7},\mathfrak{g}_{\alpha_5+\alpha_6+\alpha_7}\}.
$$ 
The subspace $Y(S_z)\otimes\mathbb{Q}$ is the orthogonal to $\epsilon_1+\epsilon_2$, $\epsilon_4-\epsilon_3$, $\epsilon_6-\epsilon_5$ and $\epsilon_7+\epsilon_8$. The restriction of $\alpha_1$ and $\alpha_3$ (resp. $\alpha_4$ and $\alpha_6$) to $Y(S_z)\otimes\mathbb{Q}$ are the two long (resp. short) simple roots of a root system of type $F_4$. Similar to the case of $E_6$, the nonzero weights for the action of $S_z$ on $\mathfrak{g}$ are precisely the roots of this root system. It follows that $C_G(H_z)^{\circ}$ is the semi-simple group of type $F_4$.
\end{proof}

\subsection{The case of type $B$ and type $C$}\label{subsec-type-BC}\hfill

Suppose that $G$ is of type $B$ or $C$. In the proof of Theorem \ref{Thm-Cent-H_z}, we have seen that there are certain root subspaces which do not contribute to the centraliser of $H_z$; namely, the root subspaces of the form (in the Lie algebras of $\Spin_{2\ell+1}$, $\Sp_{4m}$ and $\Sp_{2(2m+1)}$ respectively)
\begingroup
\allowdisplaybreaks
\begin{align}
&\{\mathfrak{g}_{\epsilon_i-\epsilon_{\ell}},\mathfrak{g}_{\epsilon_{i}},\mathfrak{g}_{\epsilon_i+\epsilon_{\ell}}\},\label{eq-subsec-type-C-G'_z-1}\\
&\{\mathfrak{g}_{2\epsilon_{2j-1}},\mathfrak{g}_{\epsilon_{2j-1}+\epsilon_{2j}},\mathfrak{g}_{2\epsilon_{2j}}\},\text{ and}\label{eq-subsec-type-C-G'_z-2}\\
&\{\mathfrak{g}_{\epsilon_{2j}-\epsilon_{\ell}},\mathfrak{g}_{\epsilon_{2j-1}-\epsilon_{\ell}},\mathfrak{g}_{\epsilon_{2j}+\epsilon_{\ell}},\mathfrak{g}_{\epsilon_{2j-1}+\epsilon_{\ell}}\}.\label{eq-subsec-type-C-G'_z-3}
\end{align}
\endgroup
Indeed, the first two sets of root subspaces realise the adjoint representation of $\SL_2$. However, if $H\subset H_z$ is an irreducible subgroup, then there is a chance that $H$ fixes a nonzero vector in the span of one of the sets of root subspaces above. 

Consider the homomorphism
$$
\Ad:\SL_2\longrightarrow\SL(\mathfrak{sl}_2).
$$
It image is isomorphic to $\PGL_2\cong\SO_3$. To describe this morphism, we fix a triple $\{x,h,y\}\subset\mathfrak{sl}_2$ satisfying the usual commuting relations and let $T_1\subset\SL_2$ be the centraliser of $h$. An element $t\in T_1\cong\mathbb{G}_m$ acts on $\mathfrak{sl}_2$ by $x\mapsto t^2x$, $h\mapsto h$ and $y\mapsto t^{-2}y$. It follows that $\Ad(T_1)$ contains regular semi-simple elements, and so $C_{\SL(\mathfrak{sl}_2)}(\Ad(T_1))$ is a maximal torus $\tilde{T}_1$.
\begin{Lem}\label{Lem-H-N_SL2(T)}
Let $H\subset\SL_2$ be an irreducible subgroup. Suppose that $\Ad(H)$ is a reducible subgroup of $\SL(\mathfrak{sl}_2)$. Then, we have $H\subset N_{\SL_2}(T_1)$ up to conjugation by $\SL_2$.
\end{Lem}
\begin{proof}
Since $H$ is an irreducible subgroup, it is linearly reductive. Suppose that $H$ has rank one, so that up to conjugation it contains $T$. Clearly, if $H=\SL_2$, then its image is not contained in any Levi subgroup. We then necessarily have $H^{\circ}=T$. It follows that $H\subset N_{\SL_2}(T)$. The only possibility is $H=N_{\SL_2}(T)$, since $T\subset\SL_2$ is not irreducible. Suppose that $H$ has rank zero; that is, a finite group. Since $H\subset\SL_2$ is irreducible, it must contain some noncentral element, say $s$, which is necessarily regular semi-simple. We may assume $s\in T$ and regard $s$ as a nonzero complex number. We first assume that $s\ne\pm\sqrt{-1}$, so that $\Ad(s)$ is regular semi-simple. By assumption, $\Ad(H)$ is a linearly reductive subgroup that is not irreducible, and thus there is a proper Levi subgroup $L$ containing $\Ad(H)$. Since $L$ contains $\Ad(s)$, it contains $\tilde{T}$. There are only four possible proper Levi subgroups containing $\tilde{T}$. If $L=\tilde{T}$, then $H$ has to be contained in $T$, contradicting the irreducibility of $H$. Now, we regard $\Ad(\SL_2)=\SO_3$ as the identity connected component of the fixed point subgroup of an involution $\theta$ acting on $\SL(\mathfrak{sl}_2)$; indeed, in the basis $\{2x,h,y\}$, the image of $\Ad$ is fixed by 
$$
\theta:g\longmapsto
\left(
\begin{array}{ccc}
& & 1\\
&1& \\
\!\!1& &
\end{array}
\right)
{}^{t}g^{-1}
\left(
\begin{array}{ccc}
& & 1\\
&1& \\
\!\!1& &
\end{array}
\right),
$$ 
where ${}^{t}g^{-1}$ is the transpose-inverse of $g$.
We denote by $\alpha_1$, $\alpha_2$ and $\alpha_1+\alpha_2$ the three positive roots of $\SL(\mathfrak{sl}_2)$ with respect to $\tilde{T}$ and the Borel subgroup $\tilde{B}$ consisting of the upper triangular matrices. The action of $\theta$ preserves $\tilde{T}$ as well as $\tilde{B}$; in particular, it permutes $\alpha_1$ and $\alpha_2$, and fixes $\alpha_1+\alpha_2$. Since $\Ad(H)\subset\SO_3$, we have $\Ad(H)\subset \theta(L)\cap L$. If $L$ has $\alpha_1$ as its unique positive roots, then $\theta(L)\cap L=\tilde{T}$. It follows that $H\subset T$, contradicting the irreducibility of $H$ again. The only possible positive root of $L$ is $\alpha_1+\alpha_2$. Consider the fixed point subgroup $L^{\theta}$, which is reductive as the centraliser of a semi-simple automorphism. Since $\theta$ acts on the root subgroup $U_{\alpha_1+\alpha_2}$ by $-1$, there is no unipotent subgroup contained in $L^{\theta}$. (Alternatively, the image of unipotent radical of a Borel of $\SL_2$ containing $T$ consists of regular unipotent elements of $\SL(\mathfrak{sl}_2)$, and thus $U_{\alpha_1+\alpha_2}$ is not contained in $\Ad(\SL_2)$.) Therefore, the identity component of $L^{\theta}$ is $(\tilde{T}^{\theta})^{\circ}$; thus, $L\cap\SO_3$ is contained in the normaliser $N_{\SO_3}((\tilde{T}^{\theta})^{\circ})$. We conclude that $H\subset N_G(T)$. Finally, we assume that $s=\pm\sqrt{-1}$ and that every other noncentral element of $H$ is conjugate to $s$; in particular, the square of any element of $H$ lies in the centre of $\SL_2$. In order for $H$ to be irreducible, there exists $s_1\in H\setminus T$. Consider the subgroup $H'\subset H$ generated by $s$ and $s_1$. An easy computation shows that $H'$ consists of the elements
$$
\{\Id,-\Id,s,-s,s_1,-s_1,ss_1,-ss_1\}.
$$ 
Obviously, the subgroup generated by $s$ is normalised by $s_1$. We deduce that $s_1\in N_{\SL_2}(T)$, since $T=C_{\SL_2}(s)$. Since this is true for any $s_1\in H\setminus T$, we have $H\subset N_{\SL_2}(T)$. Note that in this case $\Ad(H)$ is automatically reducible.
\end{proof}

\begin{Lem}\label{Lem-[t_1,t_2]}
Let $t_1$ and $t_2$ be two elements of $T_1$, and let $n\in N_{\SL_2}(T_1)$ be the matrix
$$
\left(
\begin{array}{cc}
& 1\\
-1& 
\end{array}
\right).
$$ 
Then, we have
\begin{itemize}
\item[(i)] $[t_1,t_2]=1$,
\item[(ii)] $[t_1,t_2n]=t_1^2$, $[t_1n,t_2]=t_2^{-2}$, and
\item[(iii)] $[t_1n,t_2n]=t_1^2t_2^{-2}$.
\end{itemize}
\end{Lem}
\begin{proof}
Obvious.
\end{proof}

\begin{Prop}\label{Prop-Hz-non-fix-sl2}
Let $\rho:\Pi\rightarrow\SL_2\cong H_z$ be a $z$-twisted representation. Then, there is no nonzero vector $v\in\mathfrak{sl}_2$ fixed by $\Pi$ under the action induced by $\rho$ and the adjoint representation.
\end{Prop}
\begin{proof}
Note that the centre of $\SL_2$ is the kernel of $\Ad:\SL_2\rightarrow\SL(\mathfrak{sl}_2)$, and thus a twisted $\SL_2$-representation $\rho$ defines a homomorphism $\Pi\rightarrow\SL(\mathfrak{sl}_2)$. Suppose that $v=c_1x+c_2h+c_3y\in\mathfrak{sl}_2$ is a nonzero invariant vector of $\Pi$, where $c_i$ are complex numbers, and suppose that $\rho$ is defined by the tuple $(A_i,B_i)_{1\le i\le g}$ with $A_i$ and $B_i$ lying in $N_{\SL_2}(T)\subset\SL_2$. It is not possible that $A_i$ and $B_i$ all lie in $T$, because the resulting $\rho$ is not irreducible. If $A_i=t_in$ for some $t_i\in T$ and $c_2\ne0$, where $n$ is as in Lemma \ref{Lem-[t_1,t_2]}, then $A_i$ does not fix $v$, since $t_i$ commutes with $h$ and $n$ acts on $h$ as the multiplication by $-1$. We will therefore assume that $c_2=0$. Now, $A_i$ sends $v$ to $-c_3t_i^2x-c_1t_i^{-2}y$. It follows that $c_1=-c_3t_i^2$; that is, if $A_i\notin T$, then the value of $t_i^2$ is determined by $c_1$ and $c_3$ (if they are nonzero), and the same applies to every such $B_i$. Suppose that for some $i$, both $A_i$ and $B_i$ do not lie in $T$, then Lemma \ref{Lem-[t_1,t_2]} says that $[A_i,B_i]=1$. Suppose that for some $i$, one of $A_i$ and $B_i$ does not lie in $T$, say $B_i$, and that $A_i=t_i$; then, Lemma \ref{Lem-[t_1,t_2]} says that $[A_i,B_i]=t_i^2$. However, $A_i=t_i$ fixing $v$ means that $t_i^2=1$, if $c_1$ and $c_3$ are nonzero. Therefore, $[A_i,B_i]=1$ under the assumption that $A_i\in T$ and $B_i\notin T$. Obviously, $[A_i,B_i]=1$ if both $A_i$ and $B_i$ lie in $T$. Therefore $\prod_{i=1}^g[A_i,B_i]=1$; i.e., $\rho$ is not a twisted $\SL_2$-representation. We conclude that $c_1=c_2=c_3=0$.
\end{proof}

The representation spanned by (\ref{eq-subsec-type-C-G'_z-3}) is isomorphic to $V_1\otimes V_2$, where $V_1$ (resp. $V_2$) is the standard representation of the first (resp. the second) direct factor of $H_z\cong\SL_2\times\SL_2$. 
\begin{Prop}\label{Prop-C2m+1-Stab-circ}
Suppose that $G$ is of type $C_{2m+1}$ and $z\in Z_G$ is nontrivial. Let $H'_z\cong\SL_2\subset H_z$ be the diagonal subgroup. Then, we have $C_G(H'_z)^{\circ}\cong \SO_{2m+1}$.
\end{Prop}
\begin{proof}
The maximal torus of $H_z'$ acts on the representation spanned by (\ref{eq-subsec-type-C-G'_z-3}) with weights $2$, $-2$ and $0$, and the weight space of $0$ has multiplicity two; therefore, this representation is isomorphic to the direct sum of $\mathfrak{sl}_2$ and a trivial representation, and thus contributes a root to $C_G(H'_z)^{\circ}$. Since the restriction of $\epsilon_{2j}+\epsilon_{\ell}$ to $S_z$ is a short root, we obtain a root system of type $B_m$. Note that $C_G(H'_z)^{\circ}$ and $C_G(H_z)^{\circ}$ have $S_z$ as a common maximal torus. Comparing $Y(S_z)$ and the coroot lattice of $C_G(H'_z)^{\circ}$ shows that $C_G(H'_z)^{\circ}\cong\SO_{2m+1}$. 
\end{proof}

\subsection{The connected components of $C_G(H_z)$}\label{subsec-Gamma-G_z}\hfill

For any $G$ and $z\in Z_G$, we determine the component group $C_G(H_z)/C_G(H_z)^{\circ}$ by Weyl group computations.

\begin{Lem}\label{Lem-Gamma=W/W}
Let $H\subset L_z$ be an irreducible subgroup. Then, the component group $C_G(H)/C_G(H)^{\circ}$ is isomorphic to the quotient of Weyl groups $W_{C_G(H)}/W_{C_G(H)^{\circ}}$, where $W_{C_G(H)}:=N_{C_G(H)}(S_z)/S_z$ and $W_{C_G(H)^{\circ}}:=N_{C_G(H)^{\circ}}(S_z)/S_z$.
\end{Lem}
\begin{proof}
It suffices to notice that the natural map 
$$
N_{C_G(H)}(S_z)/N_{C_G(H)^{\circ}}(S_z)\rightarrow C_G(H)/C_G(H)^{\circ}
$$ 
is surjective, which follows from the fact that every two tori of $C_G(H)^{\circ}$ are conjugate.
\end{proof}

\begin{Thm}\label{Thm-Gamma-G-z}
The component group $\Gamma:=C_G(H_z)/C_G(H_z)^{\circ}$ for each type of $G$ and $z\in Z_G$ is described as follows.
\begin{center}
\begin{tabular}{|c|c|c|c|c|c|c|c|c|c|}
\hline
$G$ & $A_{\ell}$ & $B_{\ell}$ & $C_{2m}$ & $C_{2m+1}$ & $D_{\ell}$ & $D_{2m+1}$ & $D_{2m}$ & $E_6$ & $E_7$ \\
\hline
$z$ & order $d$ & $\ne 1$ & $\ne 1$ & $\ne1$ & $\alpha_1$ & $\alpha_{\ell-1}$ or $\alpha_{\ell}$ & $\alpha_{\ell-1}$ or $\alpha_{\ell}$ & $\ne1$ & $\ne1$ \\ 
\hline
$\Gamma$ & $\mathbb{Z}/d$ & $\mathbb{Z}/2$ & $\mathbb{Z}/2$ & $\mathbb{Z}/2\times\mathbb{Z}/2$ & $\mathbb{Z}/2$ & $\mathbb{Z}/4$ & $\mathbb{Z}/2$ & $\mathbb{Z}/3$ & $\mathbb{Z}/2$. \\
\hline
\end{tabular}
\end{center}
Moreover, every element of $C_G(H_z)$ acts as an inner automorphism on $C_G(H_z)^{\circ}$ unless $G$ is of type $B$ or $C$; if $G$ is of type $B_{\ell}$ or $C_{2m}$, then there exists exactly one connected component whose elements act as an outer automorphism on $C_G(H_z)^{\circ}$; if $G$ is of type $C_{2m+1}$, then there exist two such connected components.
\end{Thm}
\begin{proof}
Suppose that $G$ is not of type $B$ or $C$. Combining Theorem \ref{Thm-Cent-H_z},  Lemma \ref{Lem-W{circ}->N(L)/L} and Remark \ref{Rem-type-Phi-proj}, we see that $W_{C_G(H_z)^{\circ}}$ is isomorphic to $N_G(L_z)/L_z$. Then, the surjective homomorphism in Lemma \ref{Lem-W->N(L)/L} shows that $W_{C_G(H_z)}$ is the semi-direct product of the normal subgroups $Z_{L_z}/S_z$ and $W_{C_G(H_z)^{\circ}}$; thus, we have $\Gamma\cong Z_{L_z}/S_z$. The connected components of $Z_{L_z}$ is determined by the torsion part of $X(T)/Q_{L_z}$, which we compute case-by-case below. 

\textit{Type $A_{\ell}$.} We may identify $Z_{L_z}$ with the subgroup of $(\mathbb{C}^{\ast})^{(\ell+1)/d}$ consisting of elements $(t_1,\ldots,t_{(\ell+1)/d})$ satisfying $\prod_{i=1}^{(\ell+1)/d}t_i^d=1$, where $S_z$ is the subtorus consisting of elements satisfying $\prod_{i=1}^{(\ell+1)/d}t_i=1$. It follows that $\Gamma$ is cyclic of order $d$.

\textit{Type $D_{\ell}$.} Suppose that $z\in Z_G$ corresponds to $\alpha_1$. The positive roots of $L_z$ are $\epsilon_{\ell-1}+\epsilon_{\ell}$ and $\epsilon_{\ell-1}-\epsilon_{\ell}$. The saturation of the root lattice of $L_z$ in $X(T)$ (i.e. the smallest sublattice of $X(T)$ containing the root lattice) contains additional vectors $\epsilon_{\ell-1}$ and $\epsilon_{\ell}$. It follows that $\Gamma\cong\mathbb{Z}/2$.

\textit{Type $D_{2m+1}$.} Suppose that $z\in Z_G$ corresponds to $\alpha_{\ell-1}$ or $\alpha_{\ell}$. The root lattice of $L_z$ is generated by $\epsilon_{2j-1}-\epsilon_{2j}$ for $1\le j\le m-1$, $\epsilon_{2m-1}-\epsilon_{2m}$, $\epsilon_{2m}-\epsilon_{2m+1}$ and $\epsilon_{2m}+\epsilon_{2m+1}$. The saturation of the lattice in $X(T)$ contains additional vectors $\epsilon_{2m}$, $\epsilon_{2m+1}$ and
$$
\frac{1}{2}(\epsilon_{2m-1}+\epsilon_{2m}+\epsilon_{2m+1}+\sum_{j=1}^{m-1}(\epsilon_{2j-1}-\epsilon_{2j})).
$$
These additional vectors generate the torsion part of $X(T)/Q_{L_z}$.

\textit{Type $D_{4m}$.} Suppose that $z\in Z_G$ corresponds to $\alpha_{\ell}$. The root lattice of of $L_z$ is spanned by $\{\epsilon_{2j-1}-\epsilon_{2j}\}_{1\le j\le 2m}$, and its saturation in $X(T)$ contains an additional vector
$$
\frac{1}{2}\sum_{j=1}^{2m}(\epsilon_{2j-1}-\epsilon_{2j}).
$$
It follows that $\Gamma\cong\mathbb{Z}/2$.

\textit{Type $E_6$.} Suppose that $z\in Z_G$ is not the identity. The root lattice of $L_z$ is spanned by $\{\alpha_1,\alpha_3,\alpha_5,\alpha_6\}$, and its saturation in $X(T)$ contains an additional vector
$$
\frac{1}{3}(\epsilon_8-\epsilon_7-\epsilon_6)-\epsilon_5=\frac{1}{3}(2\alpha_1+\alpha_3-\alpha_5-2\alpha_6).
$$
It follows that $\Gamma\cong\mathbb{Z}/3$.

\textit{Type $E_7$.} Suppose that $z\in Z_G$ is not the identity. The root lattice of $L_z$ is spanned by $\{\alpha_2,\alpha_5,\alpha_7\}$, and its saturation in $X(T)$ contains an additional vector
$$
\frac{1}{2}(\epsilon_1+\epsilon_2-\epsilon_3+\epsilon_4-\epsilon_5+\epsilon_6)=\frac{1}{2}(\alpha_2+\alpha_5+\alpha_7).
$$
It follows that $\Gamma\cong\mathbb{Z}/2$.

If $G$ is of type $B_{\ell}$ or $C_{2m}$, then $W_{C_G(H_z)^{\circ}}$ is of type $D$ according to Theorem \ref{Thm-Cent-H_z}. We have seen in Remark \ref{Rem-type-Phi-proj} that $N_G(L_z)/L_z$ is a Weyl group of type $B_{\ell-1}$ if $G$ is of type $B_{\ell}$, and is of type $B_m$ if $G$ is of type $C_{2m}$. Therefore, $W_{C_G(H_z)^{\circ}}$ is a subgroup of $N_G(L_z)/L_z$ of index two.

\textit{Type $B_{\ell}$.} The lattice $Q_{L_z}$ is spanned by $\epsilon_{\ell}$, and thus it is saturated in $X(T)$. It follows that $Z_{L_z}$ is connected. We claim that the injective homomorphism in Lemma \ref{Lem-W->N(L)/L} is surjective. Note that $N_G(L_z)/L_z$ contains the reflection $s_{\ell-1}\in W$ associated to the positive root $\epsilon_{\ell-1}$ by Lemma \ref{Lem-N_W(W_z)/W_z} and that this element does not lie in the subgroup $W_{C_G(H_z)^{\circ}}$ of type $D$. Let $n\in N_G(L_z)$ be an element lying in the connected component corresponding to $s_{\ell-1}$. Its conjugation action on $L_z$ induces an automorphism of $H_z\cong\SL_2$, which is necessarily inner, and thus there exists $n'\in N_G(L_z)$ in the same connected component as $n$ that acts on $H_z$ trivially, hence the desired surjectivity. Therefore, the Weyl group $N_{C_G(H_z)}(S_z)/S_z$ is of type $B_{\ell-1}$, and we deduce that $\Gamma\cong\mathbb{Z}/2$. The particular form of the action of $s_{\ell-1}$ on the Weyl group of $C_G(H_z)^{\circ}$ shows that the elements of $C_G(H_z)\setminus C_G(H_z)^{\circ}$ act as outer automorphisms on $C_G(H_z)^{\circ}$.

\textit{Type $C_{2m}$.} The root lattice of $L_z$ is spanned by $\{\epsilon_{2j-1}-\epsilon_{2j}\}_{1\le j\le m}$, which is saturated in $X(T)$, and thus $Z_{L_z}$ is connected. The product of reflections $w:=s_{2\epsilon_{\ell}}s_{2\epsilon_{\ell-1}}$ preserves the root system of $L_z$ and thus represents a nontrivial element of $N_G(L_z)/L_z$ by Lemma \ref{Lem-N_W(W_z)/W_z}. Let $n\in N_G(L_z)$ be an element lying the connected component corresponding to $w$. The conjugation of $n$ on $L_z$ preserves every connected component of its Dynkin diagram. Since $[L_z,L_z]$ is a product of copies of $\SL_2$, and there is no outer automorphism of $\SL_2$, there exists $n'\in N_G(L_z)$ in the same connected component as $n$ that acts on $H_z$ trivially. Since $w$ is the reflection along a coordinate vector of $Y(S_z)$, it does not lie in the type $D_m$ Weyl group $W_{C_G(H_z)^{\circ}}$. We deduce $W_{C_G(H_z)}\cong N_G(L_z)/L_z$ and the elements of $C_G(H_z)\setminus C_G(H_z)^{\circ}$ act as outer automorphisms on $C_G(H_z)^{\circ}$.

\textit{Type $C_{2m+1}$.} We first observe that the homomorphism in Lemma \ref{Lem-W->N(L)/L} is surjective with kernel $\mathbb{Z}/2$. Applying the arguments above to $w:=s_{2\epsilon_{2m}}s_{2\epsilon_{2m-1}}$ proves the surjectivity. The root lattice of $L_z$ is spanned by $\{\epsilon_{2j-1}-\epsilon_{2j}\}_{1\le j\le m}$ and $2\epsilon_{\ell}$, and thus its saturation contains an additional vector $\epsilon_{\ell}$. This shows that $Z_{L_z}/S_z\cong\mathbb{Z}/2$. We need to do slightly more to show that $\Gamma\cong\mathbb{Z}/2\times\mathbb{Z}/2$. It suffices to find an element of $W_{C_G(H_z)}$ that maps to $w$, regarded as an element of $N_G(L_z)/L_z$, but whose square belongs to $W_{C_G(H_z)^{\circ}}$. We write $H_z=H_z'\times H_z''$, where $H_z''\cong\SL_2$ has positive root $2\epsilon_{\ell}$ and $H_z'$ is the other $\SL_2$ factor. Let $L_z''$ be the Levi subgroup of $G$ containing $T$ and having $2\epsilon_{\ell}$ as its unique positive root, so that $H_z''$ is its derived subgroup. Consider $C_G(H_z'')$, and let $S_z''=Z_{L_z''}^{\circ}$, which is a maximal torus of $C_G(H_z'')$. Arguments similar to the proof of Theorem \ref{Thm-Cent-H_z} show that $C_G(H_z'')^{\circ}$ has a root system of type $C_{2m}$. Since $Y(S_z'')=\bigoplus_{i=1}^{2m}\mathbb{Z}\epsilon_i$, we have $C_G(H_z'')^{\circ}\cong\Sp_{4m}$. Let $L_z'\subset C_G(H_z'')^{\circ}$ be the Levi subgroup with positive roots $\{\epsilon_{2j-1}-\epsilon_{2j}\}_{1\le j\le m}$. Then, $[L_z',L_z']$ is a direct product of $m$-copies of $\SL_2$, and $H_z'$ is diagonally embedded in it. Write $G':=C_{C_G(H_z'')^{\circ}}(H_z')$, which contains $S_z':=Z_{L_z'}$ as a maximal torus. Our analysis of the $C_{2m}$ case above shows that $G'$ has two connected components; its Weyl group $N_{G'}(S_z')/S_z'$ is of type $B_m$ and $N_{G^{\prime\circ}}(S_z')/S_z'$ is of type $D_m$. Moreover, there exists $n'\in N_{C_G(H_z'')^{\circ}}(L_z')$ in the connected component corresponding to $s_{2\epsilon_{2m}}s_{2\epsilon_{2m-1}}$ that commutes with $H_z'$. Since $n'$ commutes with both $H_z'$ and $H_z''$, it commutes with $H_z$. However, $n'$ also normalises $L_z$ and represents $w$, and thus it represents the desired element of $W_{C_G(H_z)}$.
\end{proof}

\begin{Prop}\label{Prop-Gamma-G'-z}
Suppose that $G$ is of type $C_{2m+1}$ and $z\in Z_G$ is nontrivial. Then, we have 
$$
C_G(H'_z)/C_G(H'_z)^{\circ}\cong\mathbb{Z}/2.
$$
\end{Prop}
\begin{proof}
As in the proof of the above theorem, we deduce from Proposition \ref{Prop-C2m+1-Stab-circ} an isomorphism $W_{C_G(H_z')^{\circ}}\cong N_G(L_z)/L_z$, so that the homomorphism in Lemma \ref{Lem-W->N(L)/L} is surjective. The proof of the above theorem also shows that $Z_{L_z}/S_z\cong\mathbb{Z}/2$.
\end{proof}

\subsection{Representations with the largest stabilisers}\label{subsec-large-stab}\hfill

For any $G$ and $z\in Z_G$ nontrivial, we define $G_z:=C_G(H_z)$; if $G$ is of type $C_{2m+1}$, then we define $G'_z:=C_G(H_z')$. We will show that in each case $G_z$ is the largest stabiliser up to conjugation, and describe the representations whose stabiliser is conjugate to $G_z$ or $G_z'$.

\begin{Lem}\label{Lem-Gz-Lz}
Let $\rho\in\Rep_z(\Pi,G)$ be a completely reducible representation. Then, the maximal tori of $\Stab_G(\rho)^{\circ}$ are $G$-conjugate to a subtorus of $S_z$ and are $G$-conjugate to $S_z$ if and only if a $G$-conjugate of $\rho$ is an $L_z$-irreducible representation.
\end{Lem}
\begin{proof}
Let $S\subset \Stab_G(\rho)^{\circ}$ be a maximal torus and write $M:=C_G(S)$, and thus $\rho$ is a $z$-twisted $M$-representation. We have $L_z\subset M$ up to $G$-conjugation by Proposition \ref{Prop-L-supset-L_z}. It follows that $S\subset S_z$, and we have $S=S_z$ if and only if $M=L_z$. Note that $\rho$ is automatically $M$-irreducible by the maximality of $S$.
\end{proof}

\begin{Lem}\label{Lem-SL-L-tau}
Let $G=\SL_{\ell+1}$ and let $L\subset G$ be a Levi subgroup such that $[L,L]$ is isomorphic to the direct product of $m$-copies of $\SL_d$ with $m=(\ell+1)/d$. Let $n\in N_G(L)$ be an element that acts on $[L,L]$ by cyclically permuting the $\SL_d$-factors. Then, every $x\in Z_L^{\circ}\cap[L,L]$ can be written as $y^{-1}nyn^{-1}$ for some $y\in Z_L\cap[L,L]$.
\end{Lem}
\begin{proof}
Under the isomorphism $[L,L]\cong\SL_d\times\cdots\times\SL_d$, elements of $Z_L^{\circ}\cap[L,L]$ can be written as $x=(x_1,\ldots,x_m)$ where $x_i$'s are $d$-th roots of unity satisfying $x_1\cdots x_m=1$. For $y=(y_1,\ldots,y_m)\in[L,L]$, we may write $nyn^{-1}=(y_2,\ldots,y_m,y_1)$. The equation 
$$
x=y^{-1}nyn^{-1}=(y_2y_1^{-1},\ldots,y_my_{m-1}^{-1},y_1y_m^{-1})
$$ 
is solved by putting $y_1=1$ and solving $y_i$ inductively for all $i$. Note that $y_1y_2\cdots y_m$ is not necessarily equal to one; in other words, $y$ does not necessarily lie in $Z_L^{\circ}$.
\end{proof}

\begin{Prop}\label{Prop-Gz-Hz}
Let $\rho$ be a $z$-twisted $L_z$-representation, regarded as a point of $\Rep_z(\Pi,G)$. Suppose that $\Stab_G(\rho)$ contains $G_z$. Then, there is a $z$-twisted $H_z$-representation $\rho_1$ such that $\rho$ is a $(Z_{G_z})^{2g}$-translate of $\rho_1$.
\end{Prop} 
\begin{proof}
It follows from Theorem \ref{Thm-Cent-H_z} and the proof of Theorem \ref{Thm-Gamma-G-z} that $G_z$ meets every connected component of $N_G(L_z)$. Suppose that $\rho$ is the representation defined by $(A_i,B_i)_{1\le i\le g}$ with $A_i$ and $B_i$ lying in $L_z$. Let $n\in N_G(L_z)\cap G_z$. We have $nA_in^{-1}=A_i$ for every $i$ and similarly for the $B_i$'s. Fix some $i$ and write $A_i=zl$ with $z\in S_z$ and $l\in[L_z,L_z]$. We have $z^{-1}nzn^{-1}=lnl^{-1}n^{-1}=z_1$ for some $z_1\in S_z\cap [L_z,L_z]$. In what follows, we will choose a particular $n$ in each case and show that there exist suitable modifications of $l$ and $z$ which lie in $H_z$ and $Z_{G_z}$ respectively. Since the computation works for all $A_i$ and $B_i$ uniformly, the same conclusion holds for $\rho$.
 
\textit{Type A}. In this case, the derived subgroup $[L_z,L_z]$ is isomorphic to a direct product of $(\ell+1)/d$-copies of $\SL_d$. Let $n\in N_G(L_z)$ be as in the statement of Lemma \ref{Lem-SL-L-tau}, which commutes with $H_z$, then we obtain $y_1\in Z_{L_z}\cap [L_z,L_z]$ satisfying $y_1^{-1}ny_1n^{-1}=z_1$. We have $zl=y_1^{-1}zy_1l$, where $y_1^{-1}z\in Z_{L_z}$ and $y_1l\in[L_z,L_z]$. Since $y_1l$ commutes with $n$, it lies in $H_z$. It follows from $G_z=C_G(H_z)$ that the element $y_1l$ also commutes with $G_z$. By assumption, $zl$ commutes with $G_z$, and thus $y_1^{-1}z\in Z_{L_z}\subset C_G(H_z)$ commutes with $G_z$; therefore, $y_1^{-1}z$ lies in the centre $Z_{G_z}$.

The proof of the type-A case applies to the remaining cases. Recall the description of $H_z$ in Lemma \ref{Lem-L'-SL}. If we write $[L_z,L_z]\cong \prod_{i\in I}G_i$, then there is a subset $\tilde{I}\subset I$ and an automorphism $\sigma$ cyclically permuting the components $\{G_i\}_{i\in\tilde{I}}$, so that $H_z\subset[L_z,L_z]$ is the fixed point subgroup of $\sigma$; moreover, there is a Levi subgroup $L'\subset G$, whose derived subgroup is a special linear group, such that $\sigma$ can be realised by an element of $[L',L']$. Note that $G_i\cong\SL_{r_0}$ for every $i\in\tilde{I}$ and some $r_0$. Let $\tilde{L}\subset[L',L']$ be the Levi subgroup whose derived subgroup is $\prod_{i\in\tilde{I}}G_i$. One caveat is that Lemma \ref{Lem-SL-L-tau} assumes $x\in Z_{\tilde{L}}^{\circ}$, and thus it is not readily applicable to $z_1$, since we only know $z_1\in S_z=Z_{L_z}^{\circ}$. We must now show that the component of $z_1$ in $G_i$ for $i\notin \tilde{I}$ is trivial so that $z_1\in[L',L']$. We may choose an identification $\tilde{I}=\{1,2,\ldots,r\}$ and choose an $n\in N_G(L_z)\cap G_z$ that acts on $[L_z,L_z]$ by
$$
((g_1,g_2,\ldots,g_r),(g_i)_{i\in I\setminus\tilde{I}})\longmapsto((\sigma_r(g_r),\sigma_1(g_1),\ldots,\sigma_{r-1}(g_{r-1})),(\sigma_i(g_i))_{i\in I\setminus\tilde{I}}),
$$ 
where $\sigma_i$ is an automorphism of $\SL_{r_0}$ for every $i\in\tilde{I}$ and an automorphism of $G_i$ for every $i\in I\setminus\tilde{I}$. Since $n$ fixes every element of $H_z$, the automorphisms $\sigma_i$ must be the identity for all $i\in I$. By the equation $lnl^{-1}n^{-1}=z_1$, we deduce that the $G_i$-component of $z_1$ is trivial for every $i\in I\setminus\tilde{I}$, and that we can explicitly write $z_1$ as $(l_1l_2^{-1},l_2l_3^{-1},\ldots,l_rl_1^{-1})$, where $l_i\in G_i\cong\SL_{r_0}$ for $i\in\tilde{I}$. It follows that the product of the direct factors of $z_1$ is equal to one, and this allows us to apply the proof of Lemma \ref{Lem-SL-L-tau} to $z_1$. We thus obtain $y_1\in Z_{[\tilde{L},\tilde{L}]}\subset Z_{L_z}\cap \cap [L_z,L_z]$, which gives the desired modifications of $z$ and $l$.
\end{proof}

\begin{Cor}\label{Cor-Stab=G_z}
Suppose that $G$ is not of type $C_{2m+1}$. Let $\rho$ be a $z$-twisted $L_z$-representation, regarded as a point of $\Rep_z(\Pi,G)$. Suppose that $\Stab_G(\rho)$ contains $G_z$. Then, we have $\Stab_G(\rho)=G_z$. In particular, any $H_z$-valued $z$-twisted representation $\rho$ satisfies $\Stab_G(\rho)=G_z$.
\end{Cor}
\begin{proof}
We first show that $\Stab_G(\rho)^{\circ}=G_z^{\circ}$. Note that both $\Stab_G(\rho)^{\circ}$ and $G_z^{\circ}$ contain $S_z$ as a maximal torus. The roots of $\Stab_G(\rho)^{\circ}$ are thus a subset of $X_{(G,z)}:=\{\chi\in X(S_z)\mid\mathfrak{g}(\chi)\ne0\}$, where $\mathfrak{g}(\chi)$ is as in the proof of Theorem \ref{Thm-Cent-H_z}. It has been shown that if $G$ is not of type $B$, $C$ or $D_{2m+1}$ with $z$ corresponding to $\alpha_{\ell-1}$ or $\alpha_{\ell}$, then every $\chi\in X_{(G,z)}$ is a root of $G_z^{\circ}$; since $\Stab_G(\rho)$ contains $G_z$, they are also roots of $\Stab_G(\rho)^{\circ}$, and thus $X_{(G,z)}$ is precisely the set of roots of $\Stab_G(\rho)^{\circ}$. We deduce that $\Stab_G(\rho)^{\circ}=G_z^{\circ}$. 

Suppose that $G$ is of type $D_{2m+1}$ and $z$ corresponds to $\alpha_{\ell-1}$ or $\alpha_{\ell}$. The restrictions of the roots indexing the root subspaces (\ref{eq-Thm-Cent-H_z-D2m+1}) to $Y(S_z)$ are $\epsilon_i$. However, there are long roots $2\epsilon_i$ in the root system of $G_z$. Had the span of (\ref{eq-Thm-Cent-H_z-D2m+1}) contributed a root subspace to the Lie algebra of $\Stab_G(\rho)^{\circ}$, then its root system would be nonreduced, which is impossible. Therefore, we again have $\Stab_G(\rho)^{\circ}=G_z^{\circ}$. 

Suppose that $G$ is of type $C_{2m}$. The indexing roots in (\ref{eq-subsec-type-C-G'_z-2}) restrict to a long root defined on $Y(S_z)$. Therefore, if the Lie algebra of $\Stab_G(\rho)^{\circ}$ has a root subspace contained in the span of (\ref{eq-subsec-type-C-G'_z-2}), then it would be isomorphic to $\Sp_{2m}$. However, it is not possible for $\Sp_{2m}$ to contain $\SO_{2m}$ as a subgroup. Therefore, there are no additional roots arising from (\ref{eq-subsec-type-C-G'_z-2}) and we have $\Stab_G(\rho)^{\circ}=G_z^{\circ}$. 

Suppose that $G$ is of type $B_{\ell}$. The indexing roots in (\ref{eq-subsec-type-C-G'_z-1}) restrict to a short root defined on $Y(S_z)$. Therefore, if the Lie algebra of $\Stab_G(\rho)^{\circ}$ has a root subspace contained in the span of (\ref{eq-subsec-type-C-G'_z-2}), then it is isomorphic to $\Spin_{2\ell-1}$. Let $Z'$ be the center of $\Spin_{2\ell-1}$, which is a subgroup of $Z_{G^{\circ}_z}$, a group consisting of four elements. By Theorem \ref{Thm-Gamma-G-z}, elements of the nontrivial connected component of $G_z$ act as an outer automorphism on $G_z^{\circ}$, and thus they act nontrivially on $Z_{G_z^{\circ}}$. However, elements of $N_G(L_z)\cap G_z$ act trivially on $Z'$, since every connected component of $N_G(L_z)$ meets $\Spin_{2\ell-1}$. It follows that $Z'=Z_{G_z}$. By Proposition \ref{Prop-Gz-Hz}, there exists an $\rho_1\in\Rep_z(\Pi,H_z)$ which agrees with $\rho$ up to an element of $(Z_{G_z})^{2g}$. Since both $\rho$ and $Z_{G_z}$ commute with $\Spin_{2\ell-1}$, we  have $\Stab_G(\rho_1)^{\circ}\supset\Spin_{2\ell-1}$ (hence $\Stab_G(\rho_1)^{\circ}=\Spin_{2\ell-1}$). For any $1\le i\le \ell-1$, we regard the span of (\ref{eq-subsec-type-C-G'_z-1}) as a representation of $L_z$; for the action of the subgroup $H_z\cong\SL_2$, it is isomorphic to the adjoint representation $\mathfrak{sl}_2$. The equality $\Stab_G(\rho_1)^{\circ}=\Spin_{2\ell-1}$ means that (\ref{eq-subsec-type-C-G'_z-1}) contains a fixed vector of $\rho_1$, contradicting Proposition \ref{Prop-Hz-non-fix-sl2}. We conclude that $\Stab_G(\rho)^{\circ}=G_z^{\circ}$.

We then let $H$ be he closure of the image of $\rho$ and apply Lemma \ref{Lem-Gamma=W/W} Lemma \ref{Lem-W{circ}->N(L)/L} and Lemma \ref{Lem-W->N(L)/L} to determine the component groups; note that $\Stab_G(\rho)=C_G(H)$. Since $\Stab_G(\rho)^{\circ}=G_z^{\circ}$, the homomorphism in \ref{Lem-W{circ}->N(L)/L} is an isomorphism if $G$ is not of type $B$ or $C$. For $G$ of type $B_{\ell}$ or $C_{2m}$, the homomorphism in Lemma \ref{Lem-W->N(L)/L} is surjective. Then, Lemma \ref{Lem-W->N(L)/L} shows that $W_{C_G(H)}=W_{C_G(H_z)}$, and thus $\Stab_G(\rho)=G_z$. 
\end{proof}

\begin{Cor}\label{Cor-Stab=G_z'}
Suppose that $G$ is of type $C_{2m+1}$. Let $\rho$ be a $z$-twisted $L_z$-representation, regarded as a point of $\Rep_z(\Pi,G)$. Suppose that $\Stab_G(\rho)$ contains $G'_z$. Then, we have $\Stab_G(\rho)=G'_z$. In particular, any $H'_z$-valued $z$-twisted representation $\rho$ satisfies $\Stab_G(\rho)=G_z'$.

\end{Cor}
\begin{proof}
Both $\Stab_G(\rho)^{\circ}$ and $G_z^{\prime\circ}$ contain $S_z$ as a maximal torus. According to the proofs of Theorem \ref{Thm-Cent-H_z} and Proposition \ref{Prop-C2m+1-Stab-circ}, the only weight spaces $\mathfrak{g}(\chi)$ with $\chi\in X_{(G,z)}$ that do not contribute to the root system of $G_z'$ are the span of (\ref{eq-subsec-type-C-G'_z-2}) for some $j$. However, the corresponding weights would be long roots if they were the roots of $G'_z$, which is impossible since $G'_z$ contains short roots. We thus have $\Stab_G(\rho)^{\circ}=G_z^{\prime\circ}$. As in the proof of Corollary \ref{Cor-Stab=G_z}, we deduce $\Stab_G(\rho)=G_z^{\prime}$.
\end{proof}

For the purpose of determining the geometry of the locus where $\Stab_G(\rho)=G_z^{\prime}$ for type $C_{2m+1}$ groups, we need the following analogue of Proposition \ref{Prop-Gz-Hz}. Recall that $H_z'\subset H_z\cong\SL_2\times\SL_2$ is the diagonal subgroup. Let $H''_z\subset\SL_2\times\SL_2$ be the subgroup which has $H'_z$ as its identity connected component and has a component containing $(-\Id,\Id)$. Since $(-\Id,\Id)$ has order two and acts trivially on $H_z'$, we have $H''_z\cong H_z'\times\mathbb{Z}/2$.
\begin{Prop}\label{Prop-G'z-H'z}
Suppose that $G$ is of type $C_{2m+1}$. Let $\rho$ be a $z$-twisted $L_z$-representation, regarded as a point of $\Rep_z(\Pi,G)$. Suppose that $\Stab_G(\rho)=G'_z$. Then, $\rho$ is of the form $(\lambda_iA_i,\mu_iB_i)_{1\le i\le g}$ up to $L_z$-conjugation, where
\begin{itemize}
\item $(A_i,B_i)_i\in\Rep_z(\Pi,H''_z)$,
\item $\lambda_i$ is the nontrivial element of $Z_{G_z}$ if $A_i\notin H'_z$, and similarly for $\mu_i$ and $B_i$, and
\item $\lambda_i=1$ if $A_i\in H'_z$, and similarly for $\mu_i$ and $B_i$. 
\end{itemize}
Moreover, such $L_z$-representations have stabilisers equal to $G'_z$.
\end{Prop} 
\begin{proof}
By Proposition \ref{Prop-Gz-Hz}, $\rho$ is of the form $(\lambda_iA_i,\mu_iB_i)_{1\le i\le g}$, where $(\lambda_i,\mu_i)_i\in (Z_{G_z})^{2g}$ and $\rho_1:=(A_i,B_i)_i$ lies in $\Rep_z(\Pi,H_z)$. For any $i$, the order-two semi-simple element $\lambda_i\in S_z$ acts as $\pm 1$ on any weight space $\mathfrak{g}(\chi)$. If $\lambda_i\ne 1$, then it acts as $\eta=-1$ on $\mathfrak{g}(\chi)$ for any root $\chi$ of $G_z'$ that is not a root of $G_z$. The same is true for $\mu_i$. Consider the $\mathfrak{g}(\chi)$ spanned by (\ref{eq-subsec-type-C-G'_z-3}), so that $\chi$ is not a root of $G_z$. Let $V_1$ and $V_2$ denote the standard representation $V$ of $\SL_2$. Then, $\mathfrak{g}(\chi)$ is isomorphic to $V_1\otimes V_2$ as a representation of $H_z\cong\SL_2\times\SL_2$. Let $\{e_1,e_2\}$ be a basis of $V$. Suppose that $v_1\otimes e_1+v_2\otimes e_2$ is a nonzero vector fixed by $A_i=(g_i,h_i)$ up to $\eta$. Write
\begin{equation}\label{eq-Prop-G'z-H'z}
h_i=\left(
\begin{array}{cc}
a & b\\
c & d
\end{array}
\right)
\end{equation}
in the given basis. We have $h_ie_1=ae_1+ce_2$ and $h_ie_2=be_1+de_2$, and thus
$$
g_iv_1\otimes h_ie_1+g_iv_2\otimes h_ie_2=g_i(av_1+bv_2)\otimes e_1+g_i(cv_1+dv_2)\otimes e_2,
$$
which is required to agree with $\eta(v_1\otimes e_1+v_2\otimes e_2)$. It is not possible that $v_1$ and $v_2$ are linearly dependent. Otherwise, we may assume $v_1$ to be nonzero and $v_2=\zeta v_1$ for some $\zeta\in\mathbb{C}$, and then $g_i(a+b\zeta)v_1=\eta v_1$; however, there exists some $i$ such that $v_1$ does not span an invariant subspace of $g_i$ by the irreducibility of the representation $V$. Let $g_i'$ be the endomorphism of $V$ which acts as the matrix (\ref{eq-Prop-G'z-H'z}) in the basis $\{v_1,v_2\}$. We then have
$$
g_ig_i'^tv_1\otimes e_1+g_ig_i'^tv_2\otimes e_2=\eta (v_1\otimes e_1+v_2\otimes e_2),
$$ 
where $g_i^{\prime t}$ is the transpose, and thus $g_i=\eta g_i^{-\prime t}$. We have thus shown that $g_i$ is the transpose-inverse of $h_i$ up to $\eta$ and a base change of $V$, or rather, a conjugation by $H_z$. Moreover, the base change in question only depends on the vector $v_1\otimes e_1+v_2\otimes e_2$, and thus works uniformly for all $A_i$ and $B_i$. Now, note that the automorphism of $\SL_2$ defined by transpose-inverse is an inner automorphism, and thus a further conjugation by $H_z$ takes $(\eta h_i^{-t},h_i)$ to $(\eta h_i,h_i)$, which lies in $H''_z$ as desired.
\end{proof}
\begin{Rem}\label{Rem-G'z-H'z}
Let $z'\in S_z$ be the product of the nontrivial element of $Z_{G_z}$ and $(-\Id,\Id)\in H_z$, and let $\mathscr{Z}=\{1,z'\}$ be the group generated by $z'$. Denote by $\tilde{H}_z$ the group generated by $H'_z$ and $\mathscr{Z}$. Then, the above proposition says that the representations $\rho$ with $\Stab_G(\rho)=G'_z$ are precisely the elements of $\Rep_z(\Pi,\tilde{H}_z)$.  
\end{Rem}

\section{The most singular loci}\label{sec-sing-loci}

In this section, we give a precise meaning of the most singular locus as an algebraic variety. It is a particular Luna stratum, which we recall in \S \ref{subsec-stab-strata}. The fundamental language of \'etale slices is recall in \S \ref{subsec-EtS}. The proof of Theorem \ref{Thm-B} is given in \S \ref{subsec-geom-sing-loci}.

\subsection{\'Etale slices}\label{subsec-EtS}\hfill

We briefly review Luna's \'etale slice theorem following the lecture notes of Dr\'ezet \cite{Dre2} (see also \cite{Lu}).

Let $G$ be a reductive group acting on an affine variety $X$. Denote by $\pi:X\rightarrow X\ds G$ the affine GIT quotient. A subset $W$ of $X$ is called saturated if it is of the form $\pi^{-1}(U)$ for some subset $U\subset X\ds G$. Suppose that the orbit of $x\in X$ is closed. Then, $S:=\Stab_G(x)$ is a linearly reductive group according to Matsushima \cite{Mat} (see also  \cite[Theorem A]{Ri77}). Let $V\subset X$ be an $S$-invariant locally closed affine subvariety containing $x$. Consider the quotient variety
$$
G\times^SV:=(G\times V)\ds S,
$$
where $s\in S$ acts as $(g,v)\mapsto (gs^{-1},sv)$. Note that the action of $S$ on $G\times V$ is free. There is a natural $G$-equivariant morphism
$$
\psi:G\times^SV\longrightarrow X
$$
induced by the action of $G$ on $X$, with $G$ acting on $G\times^SV$ by left multiplication on the first factor. We say that $\psi$ is \textit{strongly \'etale} if the following conditions hold:
\begin{itemize}
\item[(i)] $\psi$ is \'etale, and its image is a saturated open subset $U\subset X$.
\item[(ii)] $\bar{\psi}:(G\times^SV)\ds G\cong V\ds S\rightarrow X\ds G$ is \'etale.
\item[(iii)] There is a Cartesian diagram
\begin{equation}\label{eq-str-et-Cart}
\begin{tikzcd}[row sep=2.5em, column sep=2em]
G\times^SV \arrow[r, "\psi"] \arrow[d, swap] & U \arrow[d, "\pi|_U"]\\
V\ds S \arrow[r, "\bar{\psi}"'] & X\ds G,
\end{tikzcd}
\end{equation}
where the left vertical arrow is the quotient by $G$.
\end{itemize}
If $\psi$ is strongly \'etale, then we say that $V$ is an \textit{\'etale slice} at $x$. Luna's theorem asserts that for any $G$, $X$ and $x\in X$ as above, an \'etale slice exists.

\subsection{Stabiliser type strata}\label{subsec-stab-strata}\hfill

A consequence of the existence of \'etale slices is the stratification of an affine GIT quotient $X\ds G$ by stabiliser types. Let $X$ be an affine algebraic variety with an action of a reductive group $G$, and let $\pi:X\rightarrow X\ds G$ be the affine GIT quotient. For any $y\in X\ds G$, there is a unique closed $G$-orbit contained in $\pi^{-1}(y)$, and let $x$ be a point in this orbit. The stabiliser $G_x=\Stab_G(x)$ is a linearly reductive group. For any closed subgroup $H\subset G$, we denote by $(H)$ the set of the $G$-conjugates of $H$, and write $(H_1)\subset(H_2)$ if a conjugate of $H_1$ is contained in $H_2$. The set $(G_x)$ only depends on $y$, and so we may write $(G_y):=(G_x)$. Define
$$
(X\ds G)_{(H)}=\{y\in X\ds G\mid (G_y)=(H)\},
$$
and
$$
(X\ds G)_{\supset(H)}=\bigcup_{(H')\supset(H)}(X\ds G)_{(H')}.
$$

These subsets have the following properties.
\begin{Prop}(\cite[Proposition 2.4]{DZ})\label{Prop-DZ}
With the notations as above, we have
\begin{itemize}
\item[(i)] for any $y\in X\ds G$, the subset $(X\ds G)_{(G_y)}$ is locally closed, and $(X\ds G)_{\supset(G_y)}$ is closed; and
\item[(ii)] there are only finitely many sets $(H)$, with $H\subset G$ being a linearly reductive subgroup, such that $(X\ds G)_{(H)}$ is nonempty.
\end{itemize}
\end{Prop}

In some specific contexts, the locally closed subvarieties $(X\ds G)_{(H)}$ are even irreducible and smooth, and satisfy $(X\ds G)_{(H_1)}\subset\overline{(X\ds G)_{(H_2)}}$ if and only if $H_1$ contains a conjugate of $H_2$. However, without any further assumption on $X$, we do not expect these properties to hold. In what follows, we will take $X=\Rep(\Pi,G)$ so that $X\ds G=\Ch(\Pi,G)$, and prove the smoothness of $(X\ds G)_{(H)}$.

\begin{Lem}\label{Lem-noncl-stab}
Let $\rho\in X$ and suppose that its $G$-orbit $G.\rho$ is not closed. Then, there exists a $\rho_0$ lying in the unique closed orbit in the closure of $G.\rho$ such that $G_{\rho}^{\circ}$ is strictly contained in $G_{\rho_0}^{\circ}$.
\end{Lem}
\begin{proof}
The closed subgroup $\Ima\rho\subset G$ admits a Levi decomposition $V\rtimes M$ where $V$ is its unipotent radical and $M$ is linearly reductive. By \cite[Proposition 2.6]{Ri88}, there exists cocharacter $\lambda$ of $G$ such that $\Ima\rho\subset P_{\lambda}$, $M\subset L_{\lambda}$ and $V\subset U_{\lambda}$. The composition of $\rho$ with the quotient $P_{\lambda}\rightarrow L_{\lambda}$ is precisely the limit of $\lambda(t)\cdot \rho$ at $t=0$, which we denote by $\rho_0$. Now $\rho_0$ lies in the closure of $G.\rho$, and its orbit is closed since $\Ima\rho_0$ is linearly reductive. By the definition of $\rho_0$, we have $G_{\rho_0}=C_G(M)$, while $G_{\rho}=C_G(\Ima\rho)$. Notice that the image of $\lambda$ is a one dimensional subgroup of $C_G(M)^{\circ}$, but it acts on $U_{\lambda}$ with nontrivial weights (and so does not commute with $V$). We see that $G_{\rho}^{\circ}$ is strictly smaller than $C_G(M)^{\circ}$.
\end{proof}

Define 
$$
X_{\supset (H)}=\{\rho\in \pi^{-1}((X\ds G)_{\supset (H)})\mid \dim G_{\rho}\ge\dim H\},
$$
and 
$$
X_{(H)}=\{\rho\in X\mid (G_{\rho})=(H)\text{ and $G.\rho$ is closed}\}.
$$
By the upper semi-continuity of stabiliser dimension, the subset $X_{\supset (H)}$ is closed in $X$. If $G.\rho$ is not a closed orbit and its closure contains a closed orbit in $X_{(H)}$, then $\rho\notin X_{\supset(H)}$ by Lemma \ref{Lem-noncl-stab}. It follows that  $$
X_{(H)}=X_{\supset (H)}\setminus\bigcup_{(H')}\pi^{-1}((X\ds G)_{\supset(H')}),
$$
where $(H')$ runs over the conjugacy classes of (linearly reductive) subgroups of $G$ that strictly contains $H$, and so $X_{(H)}$ is open in $X_{\supset(H)}$.

Define 
$$
X_{\supset H}=\{\rho\in X_{\supset(H)}\mid G_{\rho}\supset H\}
$$ 
and 
$$
X_H=\{\rho\in X\mid G_{\rho}=H\text{ and $G.\rho$ is closed}\}.
$$ 
Then, $X_{\supset H}$ is precisely the fixed point locus of $H$ in $X_{\supset(H)}$, and so is closed in $X$, while $X_H$ is open in $X_{\supset H}$ as the intersection of $X_{(H)}$ and $X_{\supset H}$.

Since $X$ is affine, the restriction of $\pi:X\rightarrow X\ds G$ to the closed subvariety $X_{\supset(H)}$ is also a GIT quotient, and the quotient is a closed subvariety of $X\ds G$. By comparing closed points, we see that $X_{\supset(H)}\ds G\cong(X\ds G)_{\supset(H)}$; therefore, $\pi$ restricts to $\pi_H:X_{\supset(H)}\rightarrow(X\ds G)_{\supset(H)}$.

\begin{Prop}\label{Prop-sm-str}
Let $H\subset G$ be a subgroup such that $(X\ds G)_{(H)}$ is nonempty. Then, $X_{(H)}$ and $(X\ds G)_{(H)}$ are smooth. Moreover, the quotient map $X_{(H)}\rightarrow(X\ds G)_{(H)}$ is a fibre bundle (\'etale locally trivial) with fibre $G/H$.
\end{Prop}
\begin{proof}
We begin by establishing the smoothness of $X_H$, from which we will deduce the smoothness of $X_{(H)}$. The final step will be showing that $X_{(H)}\rightarrow(X\ds G)_{(H)}$ is a fibre bundle with fibre $G/H$, then the smoothness of $(X\ds G)_{(H)}$ follows.

Note that $X$ is singular, so the smoothness of $X_H$ is not immediate. (In fact, $X_H$ is contained in the singular locus of $X$ as long as $H$ has positive dimension.) We will use a moduli description of $X_H$ to prove smoothness. Let $\rho\in X_H$, which represents a homomorphism $\Pi\rightarrow G$; that $\rho$ is a fixed point of $H$ is equivalent to this homomorphism factoring through $C_G(H)$. Write $\Gamma=C_G(H)/C_G(H)^{\circ}$ and let $\omega:\Pi\rightarrow\Gamma$ be the map determined by $\rho$. Denote by $\Rep_{\omega}(\Pi,C_G(H))$ the variety of homomorphisms $\Pi\rightarrow C_G(H)$ that send any element $\gamma\in\Pi$ to the prescribed connected component $\omega(\gamma)$, and it contains $\rho$ by definition. 

We deduce from the equality $G_{\rho}=H$ that $\rho$ lies in $\Rep_{\omega}^{\heartsuit}(\Pi,C_G(H))$. Denote by $C_{G}(H)^{\rho}$ the union of the connected components of $C_G(H)$ that meet $\Ima\rho$. The irreducibility of a homomorphism $\Pi\rightarrow C_G(H)^{\rho}$ means that there exists no one parameter subgroup $\lambda:\mathbb{G}_m\rightarrow C_G(H)^{\circ}$  such that $\Ima\rho$ commutes with $\Ima\lambda$ and $\Ima\lambda$ is not contained in the centre of $C_G(H)^{\rho}$; see \cite[\S 2.1 and \S 3.1]{Shu4}. Let $S\subset H^{\circ}$ be a maximal torus and let $L:=C_G(S)$, which is a Levi subgroup of $G$ containing $C_G(H)$. If $\lambda$ as above exists, then $C_L(\Ima\lambda)$ is a proper Levi subgroup $L'\subset L$ containing $\Ima\rho$. However, the connected centre $Z_{L'}^{\circ}$ is a strictly larger torus commuting with $\Ima\rho$, contradicting the maximality of $S$. We have thus proved that every element of $X_H$ lies in $\Rep_{\omega}^{\heartsuit}(\Pi,C_G(H))$ for some $\omega$. 

In view of the disjoint union
$$
\Rep(\Pi,C_G(H))=\bigsqcup_{\omega:\Pi\rightarrow\Gamma}\Rep_{\omega}(\Pi,C_G(H)),
$$
each $X_H\cap \Rep_{\omega}^{\heartsuit}(\Pi,C_G(H))$ is an open subset of $X_H$. Any $\rho\in\Rep^{\heartsuit}_{\omega}(\Pi,C_{G}(H))$ is a completely reducible $G$-representation with stabiliser containing $H$, and so $\Rep^{\heartsuit}_{\omega}(\Pi,C_{G}(H))$ is contained in $X_{\supset H}$. Since $X_H$ is open in $X_{\supset H}$, we see that $X_H\cap \Rep_{\omega}^{\heartsuit}(\Pi,C_G(H))$ is an open subset of $\Rep_{\omega}^{\heartsuit}(\Pi,C_G(H))$. Now, $\Rep^{\heartsuit}_{\omega}(\Pi,C_{G}(H))$ is smooth by \cite[Proposition 4.1]{Shu4}. It follows that $X_H$ is smooth.

Next, we prove the following equality of tangent spaces for any $\rho\in X_H$,
\begin{equation}\label{eq-Prop-sm-str}
T_{\rho}X_{(H)}=T_{\rho}X_H\oplus(\mathfrak{g}/\mathfrak{h})_H,
\end{equation}
which is similar to \cite[Lemma 3.10]{BS0}. In the above equality, $\mathfrak{g}$ and $\mathfrak{h}$ denote the Lie algebra of $G$ and $H^{\circ}$ respectively, and $(\mathfrak{g}/\mathfrak{h})_H$ is the sum of nontrivial irreducible subrepresentations in $\mathfrak{g}/\mathfrak{h}$ (since $H$ is linearly reductive, any finite dimensional representation is semi-simple). We will also use the notation $(\mathfrak{g}/\mathfrak{h})^H$ for the subrepresentation of $\mathfrak{g}/\mathfrak{h}$ on which $H$ acts trivially.

Luna's \'etale slice theorem (see \S \ref{subsec-EtS}) gives an $H$-invariant locally closed affine subvariety $V\subset X$ containing $\rho$ and a $G$-equivariant \'etale morphism $\psi:G\times^HV\rightarrow U$, where $U$ is a $G$-invariant saturated open subset of $X$. The inverse image of a closed orbit in $U$ is a finite union of closed orbits in $G\times^HV$ with conjugate stabilisers. It follows that $\psi$ restricts to an \'etale morphism $\psi_H:G/H\times V^H\rightarrow X_{(H)}\cap U$. We deduce that
$$
T_{\rho}X_{(H)}=T_{\rho}V^H\oplus(\mathfrak{g}/\mathfrak{h}).
$$
Since $X_H=X_{(H)}^H$, we have
$$
T_{\rho}X_H=(T_{\rho}X_{(H)})^H=T_{\rho}V^H\oplus(\mathfrak{g}/\mathfrak{h})^H.
$$
Then, equation (\ref{eq-Prop-sm-str}) follows:
$$
T_{\rho}X_{(H)}=T_{\rho}V^H\oplus(\mathfrak{g}/\mathfrak{h})^H\oplus(\mathfrak{g}/\mathfrak{h})_H=T_{\rho}X_H\oplus(\mathfrak{g}/\mathfrak{h})_H.
$$
In view of the fact that $X_H$ is smooth, and that every point of $X_{(H)}$ is conjugate to some point of $X_H$, the dimensions of the tangents spaces of $X_{(H)}$ are locally constant, and so $X_{(H)}$ is smooth.

It remains to show that $\pi_H:X_{(H)}\rightarrow(X\ds G)_{(H)}$ is a fibre bundle with fibre $G/H$. Applying \'etale slice to the quotient $\pi_H:X_{\supset(H)}\rightarrow(X\ds G)_{\supset(H)}$ gives the Cartesian diagram (see (\ref{eq-str-et-Cart}))
$$
\begin{tikzcd}[row sep=2.5em, column sep=2em]
G\times^HW \arrow[r, "\psi_H"] \arrow[d, swap] & X_{\supset(H)} \arrow[d, "\pi_H"]\\
W\ds H \arrow[r, "\bar{\psi}_H"'] & X_{\supset (H)}\ds G,
\end{tikzcd}
$$
where $W$ is a slice through $\rho$ in $X_{\supset(H)}$. Let $U_H$ be the image of $\psi_H$, which is a $G$-invariant saturated open subset of $X_{\supset (H)}$. Now, $U_H\cap X_{(H)}$ is a saturated open subset of $X_{\supset(H)}$, since $X_{(H)}$ is a saturated open in $X_{\supset(H)}$. Replacing $U_H$ and $W$ by respective open subsets if necessary, we have $G\times^HW\cong G/H\times W^H=G/H\times W$ and $W\ds H=W$; in other words, in an \'etale local neighbourhood of $\rho$, the map $\pi_H$ is identified with the projection $G/H\times W\rightarrow W$. The proposition is proved.
\end{proof}

\subsection{The geometry of the most singular loci}\label{subsec-geom-sing-loci}\hfill

We have seen that $G_z$ or $G'_z$ are the largest possible stabiliser groups up to conjugation. The corresponding stabiliser type strata $\Ch_z(\Pi,G)_{(G_z)}$ or $\Ch_z(\Pi,G)_{(G'_z)}$ are therefore closed subvarieties of $\Ch_z(\Pi,G)$. We will show that they are in fact a disjoint union of smooth character varieties.

\begin{Lem}
For any $z\in Z_G$, the natural inclusion $H_z\subset L_z$ induces a closed embedding
$$
\iota_{z,1}:\Ch_z(\Pi,H_z)\hooklongrightarrow\Ch_z(\Pi,L_z).
$$
Write $\bar{Z}_{G_z}:=Z_{G_z}/Z_{G_z}\cap H_z$. Denote by $\Ch_z(\Pi,Z_{G_z}\cdot H_z)\subset\Ch_z(\Pi,L_z)$ the union of the $Z_{G_z}^{2g}$-translates of $\Ch_z(\Pi,H_z)$. Then, the connected components of $\Ch_z(\Pi,Z_{G_z}\cdot H_z)$ are in bijection with $\bar{Z}_{G_z}^{2g}$.
\end{Lem}
\begin{proof}
The inclusion $[L_z,L_z]\hookrightarrow L_z$ induces a closed embedding $\Ch_z(\Pi,[L_z,L_z])\hookrightarrow\Ch_z(\Pi,L_z)$. Recall from \S \ref{subsec-Defn-H_z} that $[L_z,L_z]$ is always a direct product of special linear groups and that $H_z$ is defined as a certain diagonal subgroup. Note that both $\Rep_z(\Pi,H_z)$ and $\Rep_z(\Pi,[L_z,L_z])$ consist of stable orbits, so that two representations represent the same point of the character variety if and only if they are conjugate. Thus, it is easy to see that $\iota_{z,1}$ is an isomorphism onto the fixed point locus of some cyclic group acting on $\Ch_z(\Pi,[L_z,L_z])$. Since $Z_{G_z}$ is contained in $S_z$, there is a well-defined action of $Z_{G_z}^{2g}$ on $\Ch_z(\Pi,L_z)$. Obviously, the action of $(Z_{G_z}\cap H_z)^{2g}$ preserves the subvariety $\Ch_z(\Pi,H_z)$. Suppose that $(\lambda_i,\mu_i)_{1\le i\le g}$ is a tuple of elements  of $Z_{G_z}$ with some $\lambda_j$ not lying in $Z_{G_z}\cap H_z$. We claim that the translation of $\Ch_z(\Pi,H_z)$ by this element has empty intersection with $\Ch_z(\Pi,H_z)$, from which the second assertion follows. Nonempty intersection means that there exist two $z$-twisted $H_z$-representations $\rho=(A_i,B_i)$ and $\rho'=(A'_i,B'_i)$ and an element $l\in L_z$ such that $\lambda_jA_j=lA'_jl^{-1}$. This implies that $\lambda_j\in[L_z,L_z]$. Since $\lambda_j\in Z_{G_z}\subset S_z$, it is fixed by the conjugation action of $N_G(L_z)/L_z$; in particular, it is fixed by the cyclic permutation $\sigma$ defining $H_z$ by Lemma \ref{Lem-L'-SL}. It follows that $\lambda_j\in H_z$, which contradicts the assumption $\lambda_j\notin Z_{G_z}\cap H_z$.
\end{proof}

For $G\cong\Sp_{2(2m+1)}$ and nontrivial $z\in Z_G$, we denote by 
$$
\Ch_z(\Pi,H'_z)\hooklongrightarrow\Ch_z(\Pi,H_z)
$$
the closed subvariety induced by the inclusion $H'_z\subset H_z$, and by
$$
\Ch_z(\Pi,\mathscr{Z}\cdot H'_z)\hooklongrightarrow\Ch_z(\Pi,L_z)
$$
the union of the $\mathscr{Z}^{2g}$-translates of $\Ch_z(\Pi,H'_z)$ in $\Ch_z(\Pi,L_z)$; see Remark \ref{Rem-G'z-H'z} for the definition of $\mathscr{Z}$. Since the nontrivial element of $\mathscr{Z}$ does not lie in $H_z$, the connected components of $\Ch_z(\Pi,\mathscr{Z}\cdot H'_z)$ are in bijection with $\mathscr{Z}^{2g}$. Note that $\mathscr{Z}\cong Z_{G_z}$.
 
\begin{Thm}\label{Thm-sing-loci=Ch-A}
For any simply connected semi-simple group $G$ not of type $C_{2m+1}$ and any $z\in Z_G$, the natural inclusion $L_z\subset G$ induces a closed embedding
$$
\iota_z:\Ch_z(\Pi,Z_{G_z}\cdot H_z)\hooklongrightarrow\Ch_z(\Pi,G),
$$
and its image is precisely the closed subvariety $\Ch_z(\Pi,G)_{(G_z)}$. In particular, the connected components of $\Ch_z(\Pi,G)_{(G_z)}$ are in bijection with $\bar{Z}_{G_z}^{2g}$, and every connected component is isomorphic to $\Ch_z(\Pi,H_z)$. For $G=\Sp_{2(2m+1)}$ and the nontrivial $z\in Z_G$, the natural inclusion $L_z\subset G$ induces a closed embedding
$$
\iota'_z:\Ch_z(\Pi,\mathscr{Z}\cdot H'_z)\hooklongrightarrow\Ch_z(\Pi,G),
$$
and its image is precisely the closed subvariety $\Ch_z(\Pi,G)_{(G'_z)}$. In particular, the connected components of $\Ch_z(\Pi,G)_{(G'_z)}$ are in bijection with $Z_{G_z}^{2g}$, and every connected component is isomorphic to $\Ch_z(\Pi,H'_z)$.
\end{Thm}
\begin{proof}
By Lemma \ref{Lem-Gz-Lz}, $\Ch_z(\Pi,G)_{\supset(G_z)}$ is contained in the image of $\Ch_z(\Pi,L_z)\rightarrow\Ch_z(\Pi,G)$. We can then apply Proposition \ref{Prop-Gz-Hz} to deduce that $\Ch_z(\Pi,G)_{\supset(G_z)}$ is contained in the image of $\Ch_z(\Pi,Z_{G_z}\cdot H_z)$.  By Corollary \ref{Cor-Stab=G_z}, if $G$ is not of type $C_{2m+1}$, then every representation in the image of $\Ch_z(\Pi,Z_{G_z}\cdot H_z)$ has a stabiliser conjugate to $G_z$; thus, we see that $\Ch_z(\Pi,G)_{(G_z)}=\Ch_z(\Pi,G)_{\supset(G_z)}$ is closed and coincides with the image of $\Ch_z(\Pi,Z_{G_z}\cdot H_z)$. For $G\cong\Sp_{2(2m+1)}$, we combine Corollary \ref{Cor-Stab=G_z'} and Proposition \ref{Prop-G'z-H'z} to deduce that $\Ch_z(\Pi,G)_{(G'_z)}=\Ch_z(\Pi,G)_{\supset(G'_z)}$ and coincides with the image of $\Ch_z(\Pi,\mathscr{Z}\cdot H'_z)$.

The closed subvariety $\Ch_z(\Pi,G)_{(G_z)}$ is smooth by Proposition \ref{Prop-sm-str}. To show that $\iota_z$ is an isomorphism onto its image, it suffices to show the injectivity of $\iota_z$. Note that $\iota_z$ is a proper morphism; see, e.g., \cite[\S 3.3]{Shu4}. Indeed, a bijective proper morphism between smooth varieties is an isomorphism: for every connected component $U$ of $\Ch_z(\Pi,G)_{(G_z)}$, its inverse image $\iota_z^{-1}(U)$ is connected, since otherwise there must be a connected component of $\iota_z^{-1}(U)$ mapping surjectively onto $U$ by the properness and bijectivity of $\iota_z$, excluding the existence of other connected components by the bijectivity again; the Zariski main theorem then implies that the restriction of $\iota_z$ to $U$ is an isomorphism.

The restriction of $\iota_z$ to $\Ch_z(\Pi,H_z)$ decomposes into
$$
\Ch_z(\Pi,H_z)\stackrel{\iota_{z,1}}{\longrightarrow}\Ch_z(\Pi,L_z)\stackrel{\iota_{z,2}}{\longrightarrow}\Ch_z(\Pi,G).
$$
The second map $\iota_{z,2}$ factors through the finite quotient $\Ch_z(\Pi,L_z)/W(z)$, where $W(z):=N_G(L_z)/L_z$. Let $\rho_1$ and $\rho_2$ lie in $\Rep_z(\Pi,L_z)$. Note that they are automatically $L_z$-irreducible. Let $H_1$ (resp. $H_2$) be the Zariski closure of the image of $\rho_1$ (resp. $\rho_2$). Suppose that $g\in G$ conjugates $\rho_1$ to $\rho_2$, so that $gH_1g^{-1}=H_2$. We have $gC_G(H_1)g^{-1}=C_G(H_2)$. Note that both $C_G(H_1)^{\circ}$ and $C_G(H_2)^{\circ}$ contain $S_z=Z_{L_z}^{\circ}$ as a maximal torus. Up to modifying $g$ by an element of $C_G(H_2)^{\circ}$, we may assume that $g$ normalises $S_z$ and thus normalises $L_z$. This means that if $\iota_{z,2}(\rho_1)=\iota_{z,2}(\rho_2)$, then $\rho_1$ and $\rho_2$ represent the same point of $\Ch_z(\Pi,L_z)/W(z)$. 

To show that $\iota_z$ is injective, it suffices to show that the action of $W(z)$ on the image of $\iota_{z,1}$ is trivial; that is, for any $n\in N_G(L_z)$, and any $\rho:\Pi\rightarrow H_z$, the conjugate $n\rho n^{-1}$ can be identified with $l\rho l^{-1}$ for some $l\in L_z$. For those $G$ of types $A_{\ell}$, $D_{2m+1}$ (with $z$ corresponding to $\alpha_{\ell-1}$ or $\alpha_{\ell}$) and $E_6$, there are direct factors in $[L_z,L_z]$ that are of type $A_{\ell'}$ for some $\ell'>1$; for $G$ of other types, there are only $\SL_2$-factors. We consider the latter cases first. All automorphisms of $\SL_2$ are inner (i.e. can be realised as the conjugation by an element of $\SL_2$). It follows that the action of $n$ on $[L_z,L_z]$ is simply permuting the $\SL_2$-factors up to conjugation by some $l\in L_z$. However, $H_z$ is by definition the fixed point subgroup under a permutation that cyclically permutes all factors corresponding to short roots, and thus the action of $n$ is the same as $l$. Now, we consider other types. If $G$ is of type $A$, it is obvious that there is no element of $N_G(L_z)$ inducing an outer automorphism of a direct factor of $[L_z,L_z]$, and thus conjugation by $n$ is again a permutation combined with the conjugation by some $l\in L_z$. For the type $D_{2m+1}$, suppose that $n$ induces an outer automorphism of the $\SL_4$-factor of $[L_z,L_z]$. We may denote by $z_1$ the $\SL_4$-factor of $z$. By Proposition \ref{Prop-L-supset-L_z}, we have $z_1\ne -1$, since otherwise there exist reducible $z_1$-twisted representations. However, since any outer automorphism $\tau$ of $\SL_4$ acts as inversion on its center, we have $\tau(z_1)\ne z_1$. This contradicts the fact that $z\in Z_G$, which implies $nzn^{-1}=z$. Finally, we consider the type $E_6$. Recall $[L_z,L_z]\cong\SL_3\times\SL_3$. Let $z_1$ and $z_2$ be the two factors of $z$. Since $z$ is central, we must have $z_1=z_2$ (in fact, there is a Levi subgroup containing $L_z$ whose derived subgroup is isomorphic to $\SL_6$). Again, because an outer automorphism of $\SL_3$ acts as an inversion, there is no element of $N_G(L_z)$ inducing such an automorphism. We conclude that the action of any $n\in N_G(L_z)$ on $\rho:\Pi\rightarrow H_z$ can be identified with the conjugation of some $l\in L_z$; that is, $\iota_{z,1}$ factors through the fixed point locus of $W(z)$. Since $Z_{G_z}$ consists of fixed points of $W(z)$ in $S_z$, other connected components of $\Ch_z(\Pi,Z_{G_z}\cdot H_z)$ are also contained in the fixed point locus of $W(z)$. For $\Ch_z(\Pi,\mathscr{Z}\cdot H'_z)$, the same proof remains valid.
\end{proof}

\addtocontents{toc}{\protect\setcounter{tocdepth}{1}}
\appendix


\addtocontents{toc}{\protect\setcounter{tocdepth}{-1}}
\bibliographystyle{alpha}
\bibliography{BIB}
\end{document}